\documentclass[11pt]{article}
           
\usepackage{amsmath}
\usepackage{amsfonts}
\usepackage{amssymb}
\usepackage{mathtools}
\usepackage[mathscr]{eucal}
              \newif{\ifcomentarios}\comentariosfalse\newtheorem{theorem}{Theorem}

\newtheorem{definition}[theorem]{Definition}

\newtheorem{remark}[theorem]{Remark}

\newenvironment{proof}[1][Proof]{\textit{#1.} }{\hfill $\Box$}

\newcommand{\OBSI}{\begin{remark}\begin{rm}}
\newcommand{\OBSF}{\end{rm}\end{remark}}
\newcommand{\DEFI}{\begin{definition}\begin{rm}}
\newcommand{\DEFF}{\end{rm}\end{definition}}
\newcommand{\zerarcounters}{\setcounter{equation}{0}\setcounter{theorem}{0}}

\newcommand{\be}{\begin{eqnarray}}
\newcommand{\en}{\end{eqnarray}}
\newcommand{\bee}{\begin{eqnarray*}}
\newcommand{\ene}{\end{eqnarray*}}

\newcommand{\Z}{\mathbb{Z}}
\newcommand{\N}{\mathbb{N}}
\newcommand{\R}{\mathbb{R}}
\newcommand{\C}{\mathbb{C}}

\newenvironment{proof4}[1][Proof of Theorem~\ref{prop.espec.2}]{\noindent\textit{#1.} }{\hfill $\Box$}

\newenvironment{proof5}[1][Proof of Theorem~\ref{prop.espec.3}]{\noindent\textit{#1.} }{\hfill $\Box$}

\DeclareMathOperator{\dom}{dom}
\DeclareMathOperator{\tr}{Tr}

\DeclareMathOperator{\gerado}{span}
\newtheorem{1}{Definition}[section]

\newtheorem{3}[1]{Lemma}
\newtheorem{4}[1]{Proposition}

\newtheorem{7}[1]{Theorem}

\title{Exponential dichotomy for  dynamically defined matrix-valued Jacobi operators}
\author{Fabr\'icio Vieira Oliveira\thanks{fabricio.vieira@engenharia.ufjf.br} \and Silas Luiz Carvalho\thanks{silas@mat.ufmg.br}}
\begin{document}
	
\frenchspacing

\maketitle

% resumo em português
\begin{abstract} We present in this work a proof of the exponential dichotomy for a family of dynamically defined matrix-valued Jacobi operators in $(\mathbb{C}^{l})^{\mathbb{Z}}$, $(H_\omega)_{\omega\in\Omega}$, where $\Omega$ is a compact metric space. Namely, we show that for each $\omega\in\Omega$, the resolvent set of $H_\omega$ corresponds to the subset of $\mathbb{C}$ for which $(T,A_z)$, the $SL(2l,\mathbb{C})$-cocycle induced by the eigenvalue equation $[H_\omega u]_n=zu_n$ at $z\in\mathbb{C}$, is uniformly hyperbolic if $(T,\Omega)$ is minimal. 
\end{abstract}

\section{Introduction}
\subsection{Contextualization, definitions and main results}
%Our goal in this work is to shown the exponential dichotomy for matrix-valued Jacobi operators.
The so-called exponential dichotomy for dynamically defined self-adjoint operators consists in identifying their spectra with the subset of $\mathbb{C}$ for which the associated cocycle is non-uniformly hyperbolic; roughly, this means that  if a solution to the eigenvalue equation associated with an element of a family of dynamically defined self-adjoint operators, $(H_\omega)_\omega$ (where $\omega$ belongs to a compact metric space $\Omega$), at $z$ in the spectrum of $H_\omega$ has an exponential behavior at $\pm\infty$ (that is, it does not decay or grow exponentially at $\pm\infty$), then its rate of exponential decay or growth depends on $\omega$ (see Definition~\ref{UHC}).

Such relations between differential equations and spectral properties of the corresponding operators were firstly explored in~\cite{sell74,sell76ii,sell76iii,sell78iv}, where the concept of Sacker-Sell spectrum was presented; we also refer to \cite{palmer82} and \cite{sell80}, where the authors have explored similar relations between spectrum and hyperbolicity.

It was in~\cite{jonhson86} that these ideas were first presented for dynamically defined scalar Schr\"odinger operators, for which there is a spectral characterization of the spectrum in terms of the existence of uniform hyperbolicity for the respective cocycle. In \cite{marx14}, there is an extension of this result (by replacing the notion of uniform hyperbolicity with the notion of dominated splitting) to scalar quasiperiodic Jacobi operators (with a minimal base dynamics), including also the singular case; such result  was extended in~\cite{alkorn} for dynamically defined scalar Jacobi operators whose general base dynamics is only transitive (by extending the results of dominated splitting for $\mathbb{M}(2,\mathbb{C})$-cocycles). In \cite{haro13}, there is another extension of Johnson's theorem, now to matrix-valued operators with dynamically defined potentials, based on the results of \cite{sell76ii}.  

In the scope of the theory of dynamical systems, the equivalence between the existence of a dominated splitting and hyperbolicity is a subject of great interest; one of the most important references is \cite{shub70}. As discussed in~\cite{alkorn,marx14}, this equivalence occurs for dynamically defined scalar Jacobi operators, including the singular case. In \cite{pesin07}, a connection between uniform hyperbolicity and the existence of families of invariant cones is presented.   

There are also extensions of Johnson's theorem to CMV matrices; see~\cite{DFLY,FOZ} and references therein. 

In this work, our main goal is to show a version of the exponential dichotomy for a special class of dynamically defined matrix-valued Jacobi operators. Recall that a matrix-valued Jacobi operator in $l^{2}(\mathbb{Z}; \mathbb{C}^{l})$ is given by the law    
\begin{equation} 
\label{eq.ope.din.jacobi}
[H \textbf{u}]_{n} := D_{n - 1} \textbf{u}_{n - 1} + D_{n} \textbf{u}_{n + 1} + V_{n} \textbf{u}_{n},
\end{equation}
where $(D_{n})_{n \in \mathbb{Z}}$ and $(V_{n})_{n \in \mathbb{Z}}$ are sequences in $M(l, \mathbb{R})$ (the linear space of $l\times l$ matrices with real entries), with each $V_{n}$ being self-adjoint (see~\cite{marx15}); the case where $l=1$ corresponds to a general scalar Jacobi operator, whereas $l=1$ and $D_n=1$, for every $n\in\Z$, corresponds to a general scalar Schr\"odinger operator.

In~\cite{VieiraCarvalho2}, a characterization of the absolutely continuous spectrum (including multiplicity) of this class of matrix-like operators was presented. 

Now, let $\Omega \neq \emptyset$ be a compact metric space, $T: \Omega \rightarrow \Omega$ be a homeomorphism, % $T: \Omega \rightarrow \Omega$,
$l \in \mathbb{N}$ and let $D, V: \Omega \rightarrow M(l, \mathbb{R})$ be continuous functions. A dynamically defined matrix-valued Jacobi operator is a family of operators in $(\mathbb{C}^{l})^{\mathbb{Z}}$ given for each $\omega \in \Omega$ by the law %and each$n \in \mathbb{Z}$, by
\begin{equation} 
\label{eq.ope.din}
[H_{\omega} \textbf{u}]_{n} := D(T^{n - 1}\omega) \textbf{u}_{n - 1} + D(T^{n}\omega) \textbf{u}_{n + 1} + V(T^{n}\omega) \textbf{u}_{n}.
\end{equation}
 
If, for each $\omega\in\Omega$, $D(\omega)$ is invertible, then for each $z \in \mathbb{C}$ one may define a cocycle related with $H_\omega$ by the law
\begin{equation}
\label{eq.cociclo.din}
\begin{array}{lllll}
(T, A(z,\cdot)): & \Omega \times \mathbb{C}^{2l} & \rightarrow & \Omega \times \mathbb{C}^{2l}\\ 
&(\omega, \textbf{u}) & \mapsto & (T(\omega), A(z, \omega)\textbf{u}),
\end{array}
\end{equation}
where the map $A(z,\cdot): \Omega \rightarrow M(2l, \mathbb{C})$ is given by 
\begin{equation}
\label{def.func.mat.coci}
A(z,\omega)
:=
\left[
\begin{array}{cc}
D^{-1}(T\omega)(z - V(T\omega))    & - D^{-1}(T\omega) \\
& \\
D(T\omega) & 0
\end{array}
\right].
\end{equation}

One also has, associated with the map $A(z,\cdot): \Omega \rightarrow M(2l, \mathbb{C})$, with $z\in\mathbb{C}$, the so-called transfer matrices
\begin{equation}
\label{def.mat.trans.erg}
A_{n}(z, \omega) = 
\begin{cases}
A(z,T^{n -1}(\omega))A(z,T^{n - 2}(\omega))\cdots A(z,T(\omega))A(z,\omega), & \mbox{ if }  n \geq 1, \\
& \\
\mathbb{I}_{2l}, & \mbox{ if }  n = 0, \\
& \\
A^{-1}(z,T^{-n}(\omega))\cdots A^{-1}(z,T^{-2}(\omega))A^{-1}(z,T^{-1}(\omega)), & \mbox{ if }   n \leq - 1.
\end{cases}
\end{equation}

Each solution to the eigenvalue equation for the operator $H_{\omega}$ at $z\in\mathbb{C}$, that is, each $\mathbf{u}\in(\mathbb{C}^l)^\mathbb{Z}$ that satisfies, for each $n\in\Z$
\begin{equation}\label{eqauto}
  [H_\omega\mathbf{u}]_n=z\mathbf{u}_n,
    \end{equation}
is associated with an orbit of the cocycle; namely,  $\textbf{u}\in(\mathbb{C}^l)^\mathbb{Z}$ is a solution to the eigenvalue equation for $H_{\omega}$ at $z$ if, and only if for every $n \in \mathbb{Z}$,
\begin{equation}
\label{eq.cociclo.trans}
\left[
\begin{array}{c}
\textbf{u}_{n + 1}\\
D(T^{n}\omega) \textbf{u}_{n}
\end{array}
\right]
= 
A_{n}(z, \omega)
\left[
\begin{array}{c}
\textbf{u}_{1}\\
D(\omega) \textbf{u}_{0}
\end{array}
\right].
\end{equation}
%for every $n \in \mathbb{Z}$.
Note that for each $\omega\in\Omega$ and each $z\in\mathbb{C}$, $A(z,\omega)\in SL(2l,\mathbb{C})$, and so for each $n\in\mathbb{N}$, $A_n(z,\omega)\in SL(2l,\mathbb{C})$. % the transfer matrices~\eqref{def.mat.trans.erg} are symplectic. 

This $SL(2l,\mathbb{C})$-cocycle will be called a \textit{Jacobi cocycle}; it is an example of a more general definition. 

\begin{1}[$SL(2l,\mathbb{C})$-cocycle] \label{SLcocycle} An $SL(2l,\mathbb{C})$-cocycle over $(\Omega,\mathbb{C}^{2l})$, where $\Omega$ is a compact metric space, is defined as
\begin{equation*}
%\label{eq.cociclo.din}
\begin{array}{lllll}
(T, A): & \Omega \times \mathbb{C}^{2l} & \rightarrow & \Omega \times \mathbb{C}^{2l}\\ 
&(\omega, \textbf{u}) & \mapsto & (T(\omega), A(\omega)\textbf{u}),
\end{array}
\end{equation*}
where $T:\Omega\rightarrow \Omega$ is a homeomorphism and the map $A: \Omega \rightarrow SL(2l, \mathbb{C})$ is continuous. The iterations of $(T,A)$ are denoted by $(T,A)^n=(T^n,A_n)$, with $A_n(\omega)$ defined as in~\eqref{def.mat.trans.erg}.
\end{1}

We refer to Section~3.8 of~\cite{DF} for a complete discussion on the problem of exponential dichotomy for $SL(2,\mathbb{R})$-cocycles (and, in particular, for scalar Jacobi cocycles).

In what follows, we recall the definition of a uniformly hyperbolic $SL(2l,\mathbb{C})$-cocycle (see also~\cite{haro13}).

\begin{1}[Uniformly Hyperbolic  Cocycle]\label{UHC}
  An $SL(2l,\mathbb{C})$-cocycle  $(T,A)$ is said to be uniformly hyperbolic, and one denotes this by $(T,A)\in\mathcal{UH}$, if there exists a continuous Whitney splitting $\mathbb{C}^{2l}\times\Omega=E^s\oplus E^u$  %   such that for each $\omega\in\Omega$, $\mathbb{C}^{2l}=E^s(\omega)\oplus E^u(\omega)$, where $E^s(\omega)=s(\omega)$ and $E^u(\omega)=u(\omega)$. there exist two maps (vector bundles) $u, s : \Omega\rightarrow \mathbb{C}^{2l}$ %, $\mathbb{S}^{n-1}:=\{\textbf{u}\in\mathbb{R}^n\mid\Vert\textbf{u}\Vert=1\}$,
  such that:
\begin{enumerate}
  \item $E^s$ and $E^u$ are $(T,A)$-invariant, that is, for each $\omega\in\Omega$, $j\in\mathbb{Z}$, $\textbf{v}\in E^u(T^j\omega)$ and $\textbf{w}\in E^s(T^j\omega)$,
    \[A(T^j\omega)\textbf{v}\in E^u(T^{j+1}(\omega))\qquad \mathrm{and}\qquad  A(T^j\omega)\textbf{w}\in E^s(T^{j+1}(\omega));
    \]
\item there exists $C > 0$, $\lambda > 1$ such that for each $\omega\in\Omega$, $\textbf{v}\in E^u(\omega)$ and $\textbf{w}\in E^s(\omega)$,
  \begin{eqnarray*}
    \Vert A_{-n}(\omega)\textbf{v}\Vert\le C\lambda^{-n}\Vert\textbf{v}\Vert, \qquad \Vert A_{n}(\omega)\textbf{w}\Vert\le C\lambda^{-n}\Vert\textbf{w}\Vert, \qquad\forall\; n\in\mathbb{N}.
      \end{eqnarray*} %, for each $\textbf{v}\in u(\omega)$ and each $\textbf{w}\in s(\omega)$.
\end{enumerate}
Here, $E^u$ is called the unstable subspace and $E^s$ is called the stable subspace of $(T,A)$, and for each $\omega\in\Omega$, $E^u(\omega)$ and $E^s(\omega)$ are the fibres of $E^u$ and $E^s$, respectively, at the point $\omega$.
\end{1}

%We note that the constants $C>0$ and $\lambda>1$ do not depend on $\omega\in\Omega$, $\textbf{v}\in E^u(\omega)$ and $\textbf{w}\in E^s(\omega)$; this is why one calls such cocycles \textit{uniformly hyperbolic}. We also note that since $A:\Omega\rightarrow\mathbb{C}^{2l}$ is continuous, the fibers $E^s(\omega)$ and $E^u(\omega)$ depend continuously on $\omega\in\Omega$.

\begin{remark}\label{RUH}
  \begin{enumerate}
\item We note that the constants $C>0$ and $\lambda>1$ in Definition~\ref{UHC} do not depend on $\omega\in\Omega$, $\textbf{v}\in E^u(\omega)$ and $\textbf{w}\in E^s(\omega)$; this is why one calls such cocycles \textit{uniformly hyperbolic}. We also note that since $A:\Omega\rightarrow\mathbb{C}^{2l}$ is continuous, the fibers $E^s(\omega)$ and $E^u(\omega)$ depend continuously on $\omega\in\Omega$.%  \item It is a consequence of Definition~\ref{UHC} that there exists an invariant Whitney splitting %$\mathbb{C}^{2l}\times\Omega=E^s\oplus E^u$
\item  Let $B:\Omega\times\C\rightarrow\mathbb{C}^{2l}$ be given by the law $B(\omega,z)=A(z,\omega)$, with $A(z,\omega)$ given by \eqref{def.func.mat.coci}; then, the fibers $E^s(\omega,z)$ and $E^u(\omega,z)$ also depend continuously on $z\in\C$.
%\item Moreover, if $\Omega$ is compact, then there exists $\gamma>0$ such that  for each $n\in\Z$,
 % \[\inf_{\omega\in\Omega}\inf_{\underset{\textbf{u}(\omega)\in \textbf{U}(\omega)}{\textbf{s}(\omega)\in \textbf{S}(\omega),}}|\langle\textbf{s}(\omega),\textbf{u}(\omega)\rangle|=\inf_{\omega\in\Omega}\inf_{\underset{\textbf{u}(T^n\omega)\in \textbf{U}(T^n\omega)}{\textbf{s}(T^n\omega)\in \textbf{S}(T^n\omega),}}|\langle\textbf{s}(T^n\omega),\textbf{u}(T^n\omega)\rangle|\ge\gamma,
  %\] where for each $\omega\in\Omega$, $\textbf{S}(\omega):=s(\omega)\cap\mathbb{S}^{2l-1}${\footnote{$\mathbb{S}^{2l-1}:=\{\textbf{v}\in\mathbb{C}^{2l}\mid\Vert\textbf{v}\Vert=1\}$}} and $\textbf{U}(\omega):=u(\omega)\cap\mathbb{S}^{2l-1}$. Roughly, this means that the angle between one vector in the stable subspace and another one in the unstable subspace (both at $\omega$) is always positive. 
\end{enumerate}
\end{remark}

In order to obtain a version of the aforementioned exponential dichotomy for the family $(H_{\omega})_\omega$ given by~\eqref{eq.ope.din}, with $\Omega$ a compact metric space and $T:\Omega\rightarrow\Omega$ a minimal homeomorphism, we follow the sequence of arguments presented in \cite{zhang13} (see also Section~3.8 in~\cite{DF} for a detailed discussion of such arguments) and also explore some results that characterize minimal supports for spectral measures. Here, we consider the property known as \textit{uniform exponential growth}.

\begin{1}[Uniform Exponential Growth Condition]
\label{def.cresc.uni}
Let $(T, A)$ be an $SL(2l,\mathbb{C})$-cocycle. %with symplectic matrices $A(\omega)$ of size $2l \times 2l$.
One says that $(T,A)$ satisfies the uniform exponential growth condition, and denotes this by $(T, A) \in \mathcal{UG}$, if there exist constants $\beta > 0$ and $\lambda > 1$ such that for each fixed $\omega \in \Omega$ %, each $n\in\mathbb{N}$
and each fixed $\textbf{v}\in \mathbb{C}^{2l}\setminus\{\textbf{0}\}$,  
\begin{eqnarray}
\label{prop.cresc.exp}
\dfrac{\Vert A_{n}(\omega)\textbf{v}\Vert}{\Vert\textbf{v}\Vert}\geq \beta \lambda^{n},\;\forall\, n\in\N,\qquad \mathrm{otherwise}\qquad %s_{l}\left[A_{n}(\omega)\right] \geq \beta \lambda^{n}.
\dfrac{\Vert A_{-n}(\omega)\textbf{v}\Vert}{\Vert\textbf{v}\Vert}\geq \beta \lambda^{n},\;\forall\, n\in\N.
\end{eqnarray}
%(both conditions cannot hold simultaneously)%Here, $s_{l}\left[A_{n}(\omega)\right]$ stands for the $l$-th singular value of $A_n(\omega)$.%We denote this fact by $(T, A) \in \mathcal{UG}$
\end{1}

\begin{remark}\label{lxl}
\begin{enumerate} \item   Assume that $(T, A(z,\cdot)) \in \mathcal{UG}$ is given by~\eqref{eq.cociclo.din} and fix an orthonormal basis $\mathcal{B}:=\left\{\textbf{u}_{1}, \textbf{u}_{2},...\textbf{u}_{2l}\right\}$ of $\mathbb{C}^{2l}$ given by $2l$ linearly independent solutions %, $\{\textbf{u_1},\ldots,\textbf{u}_{2l}\}$,
    to the eigenvalue equation~\eqref{eqauto} for the operator $H_\omega$ at $z\in\mathbb{C}$. Since for each $\omega\in\Omega$, each $z\in\mathbb{C}$ and each $n\in\mathbb{N}$, $A_{n}(z,\omega)\in SL(2l,\mathbb{C})$ is a symplectic matrix, that is, $A_{n}(z,\omega)^t\mathbb{J}A_{n}(z,\omega)=\mathbb{J}$, with
    \[\mathbb{J}:=\left(\begin{array}{cc} 0&\mathbb{I}\\\mathbb{I}&0\end{array}\right),\] it follows %from the constancy of the Wronskian (Lemma~\ref{lema.const.wronsk})
      that at most $l$ vectors in $\mathcal{B}$ satisfy the condition
  \begin{equation}
\label{cond.cresc.1}
\left\|A_{n}(z,\omega)\textbf{v}\right\| \geq  \beta \lambda^{n},\qquad \forall\;n\in\N
\end{equation}
  % for each $\omega\in\Omega$ and each $n\in\mathbb{N}$;
(note that these vectors are the same for each $n\in\N$);  on the other hand, by using the same reasoning, for each $\omega\in\Omega$ %since for each $n\in\mathbb{N}$ and each $\omega\in\Omega$, $A_{-n}(T^n\omega)=[A_n(\omega)]^{-1}$,
at most $l$ vectors in $\mathcal{B}$ satisfy the condition
\begin{equation}
\label{cond.cresc.2}
\left\|A_{-n}(z,\omega)\textbf{v}\right\| \geq  \beta \lambda^{n},\qquad \forall\;n\in\N.
\end{equation}
%for each $\omega\in\Omega$ and each $n\in\mathbb{N}$.
Thus, since each vector in $\mathcal{B}$ satisfies  \eqref{cond.cresc.1} or \eqref{cond.cresc.2} (but not both simultaneously), the only possibility is that exactly $l$ vectors in $\mathcal{B}$ satisfy \eqref{cond.cresc.1} and the remaining $l$ vectors satisfy \eqref{cond.cresc.2}.
%
%Namely, it follows from the constancy of the Wronskian that for each $\omega\in\Omega$, each $z\in\mathbb{C}$ and each $n\in\mathbb{N}$, $A_{n}(z,\omega)$ is a symplectic matrix (see the proof of Lemma~\ref{lema.green.nd} for details).

\item   One has the following direct consequence of the discussion above: if $(T,A(z,\cdot))\in\mathcal{UG}$, then for each $\omega\in\Omega$, each $z\in\mathbb{C}$  and each $n\in\mathbb{Z}$, one has
\[
s_{l}\left[A_{n}(z,\omega)\right] \geq  \beta \lambda^{|n|}
\]
and 
\[s_{l + 1}\left[A_{n}(z,\omega)\right] \leq \beta^{-1} \lambda^{-|n|}, 
\]
where for each $j\in\{1,\ldots,2l\}$, $s_{j}\left[A_{n}(z,\omega)\right]$ stands for the $j$-th singular value of $A_n(z,\omega)$.
\end{enumerate}
%In fact, this reasoning applies to any basis of $\mathbb{C}^{2l}$, including any system of $2l$ linearly independent solutions %, $\{\textbf{u_1},\ldots,\textbf{u}_{2l}\}$,
%to the eigenvalue equation~\eqref{eqauto} for the operator $H_\omega$ at $z\in\mathbb{C}$.%
%This shows that for each $\omega\in\Omega$, the subspace spanned by the vectors  that satisfy~\eqref{cond.cresc.1} and the subspace spanned by the vectors  that satisfy~\eqref{cond.cresc.2} have both dimension $l$.
\end{remark}

Actually, the notion presented in Definition~\ref{def.cresc.uni} is equivalent to the notion of uniform hyperbolicity for the Jacobi cocycle $(T,A(z,\cdot))$ given by~\eqref{eq.cociclo.din} (see Subsection~\ref{UGCEUH} for a proof of this statement).

We are now able to precisely state our main result: if $\Omega$ is a compact metric space and if $T$ is a minimal homeomorphism, then for each $\omega\in\Omega$, $z\in\rho(H_{\omega})$ (where $\rho(H_\omega)$ stands for the resolvent set of $H_\omega$) if, and only if, the Jacobi cocycle $(T,A(z,\cdot))$  given by~\eqref{eq.cociclo.din} satisfies the uniform exponential growth condition. This result is a consequence of the following theorems. 

\begin{7}
\label{prop.espec.2}
Let $(H_{\omega})_{\omega}$ be the family of bounded dynamically defined operators given by~\eqref{eq.ope.din}, where $\Omega$ is a compact metric space, $T: \Omega \rightarrow \Omega$ is a minimal homeomorphism and $D,V:\Omega\rightarrow M(l,\mathbb{R})$ are continuous maps, with $D(\omega)$ invertible for each $\omega\in\Omega$. Then, for each $\omega \in \Omega$, 
\[
\rho(H_{\omega}) \subseteq \{z \in \mathbb{C}\mid (T, A(z,\cdot)) \in \mathcal{UG}\},
\]
where the associated cocycle $(T, A(z,\cdot))$ is given by \eqref{eq.cociclo.din}.
\end{7}

\begin{7}
\label{prop.espec.3}
Let $(H_{\omega})_{\omega}$ be the family of bounded dynamically defined operators given by~\eqref{eq.ope.din}, where $\Omega$ is a compact metric space, $T: \Omega \rightarrow \Omega$ is a homeomorphism and $D,V:\Omega\rightarrow M(l,\mathbb{R})$ are continuous maps, with $D(\omega)$ invertible for each $\omega\in\Omega$. Then, for each $\omega \in \Omega$, %Let $(H_{\omega})_{\omega}$ be as in the statement of Theorem~\ref{prop.espec.2}. %the family of bounded dynamically defined operators given by~\eqref{eq.ope.din}, where $\Omega$ is a compact metric space, $T: \Omega \rightarrow \Omega$ is a minimal homeomorphism and $D,V:\Omega\rightarrow M(l,\mathbb{R})$ are continuous maps, with $D(\omega)$ invertible for each $\omega\in\Omega$.
%Then, 
%for each $\omega \in \Omega$, 
\[
\rho(H_{\omega}) \supseteq \{z \in \mathbb{C}\mid (T, A(z,\cdot)) \in \mathcal{UG}\},
\]
where the associated cocycle $(T, A(z,\cdot))$ is given by \eqref{eq.cociclo.din}.
\end{7}

\begin{remark} Note that in Theorem~\ref{prop.espec.3}, there is no assumption about the dynamical nature of $(T,\Omega)$ (that it, $(T,\Omega)$ is not necessarily minimal or transitive). Theorem~\ref{prop.espec.3} corresponds to the matrix-valued version of Corollary~1 in~\cite{alkorn}.% for matrix-valued Jacobi operators.
\end{remark}
  
\subsection{Equivalence between uniform hyperbolicity and uniform growth condition for Jacobi cocycles}
\label{UGCEUH}

\begin{7}\label{UG=UH}
  Let $(T,A(z,\cdot))$ be the $SL(2l,\mathbb{C})$-cocycle given by~\eqref{eq.cociclo.din}. Then, $(T, A(z,\cdot)) \in \mathcal{UG}$ if, and only if, $(T,A(z,\cdot))\in\mathcal{UH}$. 
\end{7}
\begin{proof} One just needs to show that if $(T, A(z,\cdot)) \in \mathcal{UG}$, then $(T, A(z,\cdot)) \in \mathcal{UH}$, since the other implication follows readily from the definition of an $SL(2l,\mathbb{C})$-cocycle.

  Assume that $(T, A(z,\cdot)) \in \mathcal{UG}$ and fix an orthonormal basis $\mathcal{B}$ of $\mathbb{C}^{2l}$ given by $2l$ linearly independent solutions %, $\{\textbf{u_1},\ldots,\textbf{u}_{2l}\}$,
  to the eigenvalue equation~\eqref{eqauto} for the operator $H_\omega$ at $z\in\mathbb{C}$. Now, for each $\omega\in\Omega$, let $E^s(\omega):=\gerado(\mathcal{B}\setminus\mathcal{B}_\omega)$, where $\mathcal{B}_\omega$ are the elements of $\mathcal{B}$ that satisfy~\eqref{cond.cresc.1} (note that, by Remark~\ref{lxl}, $\mathcal{B}_\omega\notin\{\emptyset,\mathcal{B}\}$; in fact, $\#(\mathcal{B}_\omega)=l$). Accordingly, for each $\omega\in\Omega$, let $E^u(\omega):=\gerado(\mathcal{B}_\omega)$. %, where $\mathcal{C}_\omega$ are the elements of $\mathcal{B}$ that satisfy~\eqref{cond.cresc.2}.

  From now on, we omit the dependence on $z$. Let us show that $E^u$ and $E^s$ satisfy the conditions stated in Definition~\ref{UHC}.  Namely, since for each $\omega\in\Omega$ and each $n\in\mathbb{N}$, $A_{n}(\omega)\in SL(2l,\mathbb{C})$, it follows from~\eqref{prop.cresc.exp},~\eqref{cond.cresc.1} and~\eqref{cond.cresc.2} that for each $\omega\in\Omega$, $\textbf{v}\in E^s(\omega)\setminus\{\textbf{0}\}$ and $\textbf{w}\in E^u(\omega)\setminus\{\textbf{0}\}$,
\begin{equation}\label{cond.hip.uni.}
  \dfrac{\Vert A_{-n}(\omega)\textbf{w}\Vert}{\Vert\textbf{w}\Vert},\dfrac{\Vert A_{n}(\omega)\textbf{v}\Vert}{\Vert\textbf{v}\Vert}\le C\lambda^{-n},\qquad\forall\; n\in\N,
  \end{equation}
where $C:=c\beta^{-1}$, with $c$ a positive constant such that for each $\textbf{u}\in\mathbb{C}^{2l}$,  $\Vert \textbf{u}\Vert_1\le c\Vert\textbf{u}\Vert_2$ (such constant exists, since the norms $\Vert\textbf{u}\Vert_1:=\max\{|u_1|,\ldots,|u_{2l}|\}$ and $\Vert\textbf{u}\Vert_2:=(\sum_j|u_j|^2)^{1/2}$ are equivalent).
This proves condition~2.

Now, suppose that there exist $\omega\in\Omega$, $j\in\mathbb{Z}$ and $\mathbf{0}\neq\textbf{v}\in E^s(T^j\omega)$ so that $A(T^j\omega)\textbf{v}\notin E^s(T^{j+1}\omega)$; given that $\mathbb{C}^{2l}=E^s(T^{j+1}\omega)\oplus E^u(T^{j+1}\omega)$, one has $A(T^j\omega)\textbf{v}=a\textbf{w}+b\textbf{u}$ for  some $\textbf{w}\in E^u(T^{j+1}\omega)$ and $\textbf{u}\in E^s(T^{j+1}\omega)$, with $a\neq 0$. %)ausing the reasoning above, there exists $\textbf{w}\in u(T\omega)$ such that $(\textbf{w},A(\omega)\textbf{v})\neq 0$. Let $A(\omega)\textbf{v}$.
It follows from the definition of $E^s(T^{j+1}\omega)$ and $E^u(T^{j+1}\omega)$ that for each $n\in\mathbb{N}$,
\[\dfrac{\Vert A_n(T^{j+1}\omega)A(T^j\omega)\textbf{v}\Vert}{\Vert A(T^j\omega)\textbf{v}\Vert}=\dfrac{\Vert A_n(T^{j+1}\omega)(a\textbf{w}+b\textbf{u})\Vert}{\Vert A(T^j\omega)\textbf{v}\Vert}\ge c_1\lambda^n-c_2\lambda^{-n},\]
with $c_1:=C^{-1}(\Vert a\textbf{w}\Vert/\Vert A(T^j\omega)\textbf{v}\Vert)>0$, $c_2:=C(\Vert b\textbf{u}\Vert/\Vert A(T^j\omega)\textbf{v}\Vert)\ge0$. On the other hand, it follows from~\eqref{cond.hip.uni.} that for each $n\in\mathbb{N}$,
\[\Vert A_{n}(T^{j+1}\omega)A(T^j\omega)\textbf{v}\Vert=\Vert A_{n+1}(T^j\omega)\textbf{v}\Vert\le  C\lambda^{-n-1}\Vert\textbf{v}\Vert;\]
thus, for each $n\in\mathbb{N}$,
\[c_1\le \lambda^{-2n}(c_2+c_3),\]
where $c_3:=C\lambda^{-1}(\Vert\textbf{v}\Vert/\Vert A(T^j\omega)\textbf{v}\Vert)>0$. % which is an absurd.
From this contradiction, one concludes that for each $\omega\in\Omega$, $j\in\mathbb{Z}$ and $\textbf{v}\in E^s(T^j\omega)$, $A(T^j\omega)\textbf{v}\in E^s(T^{j+1}\omega)$.

Using the same reasoning, one may prove that for each $\omega\in\Omega$, $j\in\mathbb{Z}$ and and $\textbf{w}\in E^u(T^j\omega)$, $A(T^j\omega)\textbf{w}\in E^u(T^{j+1}\omega)$. This shows condition~1, and we are done.
\end{proof}  

\begin{remark} Actually, the result stated in Theorem~\ref{UG=UH} is also valid for more general $SL(2l,\mathbb{C})$-cocycles. Namely, let $\mathcal{C}_\omega$ be the set of vectors $\textbf{v}\in\mathbb{C}^{2l}$  that satisfy
\begin{equation*}
%\label{cond.cresc.1}
\frac{\left\|A_{n}(z,\omega)\textbf{v}\right\|}{\left\|\textbf{v}\right\|} \geq  \beta \lambda^{n},\qquad \forall\;n\in\N
\end{equation*}
(which is non-empty, given that the cocycle is $SL(2l,\mathbb{C})$),  let $\mathcal{B}_\omega^u$ be a maximal linearly independent subset of $\mathcal{C}_\omega$ and let, for each $\omega\in\Omega$, $E^u(\omega):=\gerado(\mathcal{B}_\omega^u)$.

Accordingly, let $\mathcal{B}_\omega^s$ be a maximal linearly independent subset of $\mathbb{C}^{2l}\setminus \mathcal{C}_\omega$ and let, for each $\omega\in\Omega$, $E^s(\omega):=\gerado(\mathcal{B}_\omega^s)$. Naturally,  it follows from Definition~\ref{def.cresc.uni} that  for each $\omega\in\Omega$, $\mathbb{C}^{2l}=E^s(\omega)\oplus E^u(\omega)$.

The rest of the proof follows the same steps presented in the proof of Theorem~\ref{UG=UH}.
  \end{remark}

The exponential dichotomy for dynamically defined matrix-valued Jacobi operators is now a direct consequence of
Theorems~\ref{prop.espec.2},~\ref{prop.espec.3} and~\ref{UG=UH}.

\begin{7}\label{maintheo}
Let $(H_{\omega})_{\omega}$ be as in the statement of Theorem~\ref{prop.espec.2}. Then, %the family of bounded dynamically defined operators given by~\eqref{eq.ope.din}, where $\Omega$ is a compact metric space, $T: \Omega \rightarrow \Omega$ is a minimal homeomorphism and $D,V:\Omega\rightarrow M(l,\mathbb{R})$ are continuous maps with $D(\omega)$ invertible for each $\omega\in\Omega$. Then, %with $\dom(H_\omega)=l^{2}(\mathbb{Z}, \mathbb{C}^{l})$ for each $\omega\in\Omega$. If the associated cocycle $(T, A)$, given by \eqref{eq.cociclo.din}, is defined in a compact metric space $\Omega$ and if $T: \Omega \rightarrow \Omega$ is a minimal homeomorphism, then,
for each $\omega \in \Omega$, 
$$
\rho(H_{\omega})=\{z \in \mathbb{C}\mid (T, A(z,\cdot)) \in \mathcal{UG}\}=\{z \in \mathbb{C}\mid (T, A(z,\cdot)) \in \mathcal{UH}\},
$$
with the associated cocycle $(T, A(z,\cdot))$  given by \eqref{eq.cociclo.din}. 
\end{7}

\begin{remark} Given that $(H_\omega)_\omega$ is a family of self-adjoint operators, it follows that for each $\omega\in\Omega$, $\sigma(H_\omega)$ (the spectrum of $H_\omega$) is a subset of $\mathbb{R}$. Therefore, the result stated in Corollary~\ref{maintheo} is equivalent to the statement that for each $\omega\in\Omega$,\[\sigma(H_{\omega})=\{x \in \mathbb{R}\mid (T, A(x,\cdot)) \notin \mathcal{UG}\}.
  \]
\end{remark}

We organize this paper as follows. In Section~\ref{crescimento} we present a necessary condition for an $SL(2l,\mathbb{C})$-cocycle to satisfy the uniform exponential growth condition (this is Proposition~\ref{prop.caracter}). We also show that the uniform exponential growth condition for the cocycle $(T,A(z,\cdot))$ is open with respect to $z\in\C$ (this is Proposition~\ref{prop.hip.aberto}).

In Section~\ref{expdic} we present the proofs of Theorems %we split the proof of Theorem~\ref{maintheo} in two parts; firstly, we prove for each $\omega\in\Omega$ the set inclusion $\rho(H_{\omega})\subseteq\{z \in \mathbb{C}\mid (T, A_z) \in \mathcal{UG}\}$ (this is Theorem~
\ref{prop.espec.2} and %), and then, we prove for each $\omega\in\Omega$ the set inclusion $\rho(H_{\omega})\supseteq\{z \in \mathbb{C}\mid (T, A_z) \in \mathcal{UG}\}$ (this is Theorem~
\ref{prop.espec.3}.

\section{$SL(2l,\mathbb{C})$-cocycles and the Uniform Exponential Growth Condition}
\label{crescimento} 
\zerarcounters

In this section, we present a necessary condition for an $SL(2l,\mathbb{C})$-cocycle $(T,A)$ (see Definition~\ref{SLcocycle}) %defined over a compact metric space $\Omega$, with $A: \Omega \rightarrow SL(2l, \mathbb{C})$ a continuous map,
to satisfy $(T,A)\notin\mathcal{UG}$ (Proposition \ref{prop.caracter}). This result is required in the proof of Theorem~\ref{prop.espec.2}. We adapt some of the ideas employed in \cite{zhang13} for $SL(2, \mathbb{R})$-cocycles (see also Section~3.8 in~\cite{DF}).

We also show that the condition $(T,A(z,\cdot))\in\mathcal{UG}$ in open for $z\in\C$ (Proposition~\ref{prop.hip.aberto}), where $(T,A(z,\cdot))$ stands for the cocycle defined by~\eqref{def.mat.trans.erg}. This result is required in the proof of Theorem~\ref{prop.espec.3}.

\begin{3}
\label{lema.carac.omega.v}
Let $(T, A)$ be a cocycle defined over a compact metric space $\Omega$, with $A: \Omega \rightarrow SL(2l, \mathbb{C})$ a continuous map. If there exist $\epsilon > 0$ and $R \in \mathbb{N}$ such that for every $(\omega, \textbf{v}) \in \Omega \times \mathbb{S}^{2l-1}$, there exists $r \in \mathbb{Z}$ such that $\left|r\right| \leq R$ and
\begin{equation}
\label{eq.matr.trans.l}
\left\| A_{r}(\omega)\textbf{v} \right\| \geq 1 + \epsilon,
\end{equation}
then $(T, A)\in\mathcal{UG}$. % satisfies the uniform exponential growth condition.
\end{3}
\begin{proof}
Let $(\omega, \textbf{v}) \in \Omega \times \mathbb{S}^{2l-1}$ and set %Suppose that $(\omega, \textbf{v})$ satisfies \eqref{eq.matr.trans.l}. Let, for each fixed denote by $r(\omega, \textbf{v})$ the number given by
$$
\begin{array}{lll}
r(\omega, \textbf{v}) & := & m; \qquad\left|m\right| =  \min \{\left| r \right| \mid \left| r \right| \leq R, $ where $ r, \omega $ and $ \textbf{v} $ satisfy \eqref{eq.matr.trans.l}$ \}
\end{array}
$$
(by hypothesis, there exists at least one $r\in\mathbb{Z}$, with $|r|\le R$, for which \eqref{eq.matr.trans.l} is valid for each $(\omega, \textbf{v}) \in \Omega \times \mathbb{S}^{2l-1}$; moreover, $r(\omega, \textbf{v})\neq 0$).

One defines recursively a sequence $(r_{k}, \textbf{v}_{k}, \omega_{k})_{k}$ in $\mathbb{Z} \times \mathbb{S}^{2l-1} \times \Omega$ by the following procedure: for $k = 0$, let $r_{0} = 0$, $\textbf{v}_{0} = \textbf{v}$, and $\omega_{0} = \omega$; for $k > 0$, set
$$ 
\left\{
\begin{array}{lll}
r_{k} & := & r(\omega_{k - 1}, \textbf{v}_{k - 1}), \\
& & \\
\textbf{v}_{k} & := & \dfrac{A_{r_{k}}(\omega_{k - 1})\textbf{v}_{k - 1}}{\left\| A_{r_{k}}(\omega_{k - 1})\textbf{v}_{k - 1} \right\|}, \\
& & \\
\omega_{k}& := & T^{r_{k}}(\omega_{k - 1}).
\end{array}\right.
$$

Thus, for every $k \in \mathbb{N}$,
$$
\left\| A_{r_{k  + 1}}(\omega_{k}) \textbf{v}_{k} \right\| \geq 1 + \epsilon.
$$

Now, note that by the definition of $\textbf{v}_{k}$, 
$$
\left\|A_{r_{k + 1}}(\omega_{k}) \frac{A_{r_{k}}(\omega_{k - 1})\textbf{v}_{k - 1}}{\left\| A_{r_{k}}(\omega_{k - 1})\textbf{v}_{k - 1} \right\|}\right\| \geq 1 + \epsilon;
$$
it follows from~\eqref{def.mat.trans.erg} that
$$
A_{r_{k + 1}}(\omega_{k}) A_{r_{k}}(\omega_{k - 1})\textbf{v}_{k - 1}
= 
A_{r_{k + 1} + r_{k}}(\omega_{k - 1})\textbf{v}_{k - 1},
$$
and so
%$$
%\frac{\left\|A_{r_{k + 1} + r_{k}}(\omega_{k - 1})\textbf{v}_{k - 1}\right\|}{\left\| A_{r_{k}}(\omega_{k - 1})\textbf{v}_{k - 1} \right\|} \geq 1 + \epsilon,
%$$
%that is, 
$$
\left\|A_{r_{k + 1} + r_{k}}(\omega_{k - 1})\textbf{v}_{k - 1}\right\| \geq (1 + \epsilon)\left\| A_{r_{k}}(\omega_{k - 1})\textbf{v}_{k - 1} \right\|\ge(1 + \epsilon)^{2}.
$$

If one proceeds recursively, one gets, for $0\le j\le k-1$,
\begin{equation}
\label{eq.rec.lk}
\left\|A_{r_{k + 1} + \ldots + r_{k - j}}(\omega_{k - j - 1})\textbf{v}_{k - j - 1}\right\| \geq (1 + \epsilon)^{j + 2}.
\end{equation}

Set, for each $k \in \mathbb{Z}_0^+$, $R_{k} := \sum^{k}_{j = 0}r_{j}$; it follows from relation \eqref{eq.rec.lk} that for each $k\in\mathbb{N}$ and each $1\le j\le k$,
$$
\left\| A_{R_{k} - R_{j - 1}}(\omega_{j - 1})\textbf{v}_{j - 1} \right\| \geq (1 + \epsilon)^{k - j+1},
$$
and in particular for $j=1$, one gets
\begin{equation}\label{crescente}
\left\| A_{R_{k}}(\omega)\textbf{v} \right\| \geq (1 + \epsilon)^{k}.
\end{equation}

It follows from~\eqref{crescente} that there exists $L$ such that for $k > L$, $\left|R_{k}\right| > R$. Indeed, let $M:=\max_{\omega\in\Omega}\max_{|l|\le R}\Vert A_l(\omega)\Vert<\infty$ and set $L:=\max\{k\in\mathbb{N}\mid (1+\epsilon)^{k-1}\le M\}$; then, the assertion follows. %one has from~\eqref{crescente} that for each $(\omega,\textbf{v})\in \Omega\times\mathbb{S}^{2l-1}$ and each $k\ge L$, $|R_k|>R$. %there exists $L=L(\omega,\textbf{v})\in\mathbb{N}$ such that for each $k>L$, $\left\| A_{R_{k}}(\omega)\textbf{v}\right\|>M$. Let $L_0:=\max\{|k|\le L\mid |R_k|\le R\}$, and note that since the map $\Omega\times \mathbb{S}^{2l-1}\ni (\omega,\textbf{v})\mapsto \Vert A_n(\omega)\textbf{v}\Vert\in\mathbb{R}$ is continuous,
%
%Set $k_0:=\max_{\omega\in\Omega}\max_{\textbf{v}\in\mathbb{S}^{2l-1}}L_0(\omega,\textbf{v})$

At last, set $k_0:=\max\{k\le L\mid |R_k|\le R\}$. Note that by relation~\eqref{crescente}, $|R_k|\to\infty$; moreover, by the definition of $L$, one has that $R_k>R$ for each $k>k_0$, otherwise $R_k<-R$ for each $k>k_0$. By combining these facts with $\left|r_{k}\right|\le R$, valid for each $k\in\mathbb{N}$, one concludes that there exists a sequence $\{k_j\}$, with $k_1 > k_{0}$, such that $R_{k_{j + 1}} > R_{k_j}$, $R_{k_{j + 1}}- R_{k_j}\le R$, otherwise there exists a sequence $\{l_j\}$, with $l_1>k_0$, such that $R_{l_{j + 1}} < R_{l_j}$, $R_{l_{j}}- R_{l_{j+1}}\le R$.

Suppose initially that $R_{k_j + 1} > R_{k_j}$ for every %$k > k_{0}$.
$j\in\mathbb{N}$. Then,  %there exists $j_0\in\mathbb{N}$ such that
$R_{k_1}>R\ge R_{k_0}>0$, and for each $n > R_{k_0}$, there exists $Q \in \mathbb{Z}_0^+$ such that $R_{k_Q} \leq n < R_{k_{Q + 1}}$; thus, one can write
\[\Vert A_{R_{k_Q}-n}(T^{n}\omega)\Vert\left\|A_{n}(\omega)\textbf{v}\right\| 
 \geq \left\| (A_{n - R_{k_Q}} (T^{R_{k_Q}}\omega))^{-1}   A_{n}(\omega)\textbf{v}\right\| 
=% \geq   
%s_{2l}\left[A_{n - R_{Q}} (T^{R_{Q}}\omega)\right]
\left\|A_{R_{k_Q}}(\omega)\textbf{v}\right\|.
\]
%where it was used the fact that, for each $A,B\in M(2l,\mathbb{C})$, $\Vert AB\Vert\ge s_{2l}[A]s_1[B]=s_{2l}[A]\Vert B\Vert$ (see~\cite{wang97}). 

Now, given that %since $n - R_{Q} \leq R$, %it follows from % and matrices are symplectic, by the
%the fact that
the map $A:\Omega\rightarrow SL(2l,\mathbb{C})$ is continuous (with $A(\omega)$ invertible for each $\omega\in\Omega$), $\Omega$ is compact and $n - R_{k_Q} <R_{k_{Q+1}}-R_{k_Q}\leq R$, there exists a constant $D > 0$ such that for each $n>R_{k_0}$ and each $Q\in\mathbb{N}$ so that $n-R_{k_Q}\leq R$, one has  
$$
\Vert A_{R_{k_Q}-n} (T^{n}\omega)\Vert\le(\max_{\omega\in\Omega}\Vert A(\omega)^{-1}\Vert)^R=:D.
$$
Therefore, %taking into account that for each $j\in\mathbb{N}$, $0<r_{k_j} \leq R$,
it follows from the estimates
\[n-R<R_{k_{Q+1}}-R\le R_{k_{Q+1}}-R_{k_0}\le R_{k_Q}+R-R_{k_0}=\sum_{k=k_0}^{k_{Q}}r_{k}+R\le(k_{Q}-k_0+1)R
\]
that $k_{Q} + 2-k_0 > \frac{n}{R}$, and so
$$
\left\|A_{n}(\omega)\textbf{v}\right\| 
 \geq   
G \left\|A_{R_{k_Q}}(\omega)\textbf{v}\right\|
\geq
G(1 + \epsilon)^{k_{Q}}
\geq
G(1 + \epsilon)^{k_{0} + \frac{n}{R} - 2}\ge G(1 + \epsilon)^{\frac{n}{R} - 2},
$$
with $G=D^{-1}$.

Suppose now that $R_{l_{j + 1}} < R_{l_j}$ for every $l_1 > k_{0}$. In this case, $R_{l_1}<-R\le R_{k_0}<0$, and for each $n\in\mathbb{N}$ such that $-n < R_{k_0}$, there exists $Q \in \mathbb{Z}_0^+$ such that  $R_{l_{Q + 1}} < -n \leq R_{l_Q}$. Then, it follows from arguments analogous to those used in the previous case that there exists $\tilde{G}>0$ such that, for each $n\in\mathbb{N}$ so that $-n < R_{k_{0}}$, one has
$$
\left\|A_{-n}(\omega)\textbf{v}\right\| 
 \geq   
%C \left\|A_{R_{Q}}(\omega)\textbf{v}\right\|
%\geq
%C(1 + \epsilon)^{Q - 1}
%\geq
\tilde{G}(1 + \epsilon)^{k_{0} + \frac{n}{R} - 2}\ge\tilde{G}(1 + \epsilon)^{\frac{n}{R} - 2}.
$$

Now, set $\lambda:=(1+\epsilon)^{1/R}$ and $C:=\min\{m/\lambda^R,\min\{G,\tilde{G}\}(1+\epsilon)^{-2}\}$, where $m:=\min_{\omega\in\Omega}(\Vert A(\omega)\Vert^{-R})>0$; %\min_{|k|\le R}\Vert A_k(\omega)\textbf{v}\Vert>0$;
then, %Summing up, there exist constants $C>0$ and $\lambda>1$ such that
it was shown that for each fixed $\omega \in \Omega$ and $\textbf{v} \in \mathbb{S}^{2l-1}$, %  and each $n\in\mathbb{N}$, 
\begin{equation*}
%\label{cond.cresc.1}
\left\|A_{n}(\omega)\textbf{v}\right\| \geq  C \lambda^{n},\qquad \forall\,n\in\mathbb{N},
\end{equation*}
otherwise
\begin{equation*}
%\label{cond.cresc.2}
\left\|A_{-n}(\omega)\textbf{v}\right\| \geq  C \lambda^{n},\qquad \forall\,n\in\mathbb{N}.
\end{equation*}
\end{proof}

%\begin{remark} In fac
%Now fix an orthonormal basis $\mathcal{B}:=\left\{\textbf{u}_{1}, \textbf{u}_{2},...\textbf{u}_{2l}\right\}$ of $\mathbb{C}^{2l}$. Since for each $\omega\in\Omega$ and each $n\in\mathbb{N}$, $A_{n}(\omega)\in SL(2l,\mathbb{R})$, at most $l$ vectors in $\mathcal{B}$ satisfy the condition \eqref{cond.cresc.1}. On the other hand, using the same reasoning, %since for each $n\in\mathbb{N}$ and each $\omega\in\Omega$, $A_{-n}(T^n\omega)=[A_n(\omega)]^{-1}$,
%at most $l$ vectors in $\mathcal{B}$ satisfy the condition \eqref{cond.cresc.2}. Thus, since each vector in $\mathcal{B}$ satisfies  \eqref{cond.cresc.1} or \eqref{cond.cresc.2}, the only possibility is that exactly $l$ vectors satisfy \eqref{cond.cresc.1} and the remaining $l$ vectors satisfy \eqref{cond.cresc.2}. Therefore, for each $\omega\in\Omega$ and each $n\in\mathbb{N}$, one has
%$$
%s_{l}\left[A_{n}(\omega)\right] \geq  C \lambda^{n}.
%$$
%\end{proof}

\begin{4}
\label{prop.caracter}
Let $(T, A)$ be a cocycle defined over a compact metric space $\Omega$, with $A: \Omega \rightarrow SL(2l, \mathbb{C})$ a continuous map. If $(T, A) \notin \mathcal{UG}$, then there exist $(\omega,\textbf{v})\in \Omega\times\mathbb{S}^{2l-1}$ such that for each $n\in\mathbb{Z}$, 
$$
\left\|A_{n}(\omega)\textbf{v}\right\| \le 1.
$$
\end{4}
\begin{proof}
If  $(T, A)\notin \mathcal{UG}$, it follows from Lemma \ref{lema.carac.omega.v} that for each $k \in \mathbb{N}$, there exists a pair $(\omega_{k}, \textbf{v}_{k}) \in \Omega \times \mathbb{S}^{2l-1}$ such that for each $r \in \mathbb{Z}$ with $\left| r \right| \leq k$, one has
$$
\left\|A_{r}(\omega_{k})\textbf{v}_{k}\right\| < 1 + \frac{1}{k}.
$$ 
Since $\Omega\times \mathbb{S}^{2l-1}$ is compact, there exist $(\omega,\textbf{v}) \in \Omega\times\mathbb{S}^{2l-1}$ and sequences $(\omega_{k_{j}})_{j}$ and $(\textbf{v}_{k_{j}})_{j}$ such that 
$$
\begin{array}{lll}
\lim_{j \rightarrow \infty} \omega_{k_{j}} & = & \omega, \\
&& \\
\lim_{j \rightarrow \infty} \textbf{v}_{k_{j}} & = & \textbf{v}. 
\end{array}
$$
It follows now from the fact that the map $A:\Omega\rightarrow SL(2l,\mathbb{C})$ is continuous that for each $n\in\mathbb{Z}$,
$$
\left\|A_{n}(\omega)\textbf{v}\right\| \leq \lim_{k\to\infty}\Vert A_n(\omega_k)\textbf{v}_k\Vert\le \lim_{k\to\infty}\left(1+\dfrac{1}{k}\right)=1.
$$
%
%Now, using the same reasoning as before, there exist $\omega \in \Omega$ and $\textbf{v} \in S(\textbf{0};1)$ such that 
%$$
%\begin{array}{lll}
%\lim_{j \rightarrow \infty} \omega_{R_{j}} & = & \omega, \\
%&& \\
%\lim_{j \rightarrow \infty} \textbf{v}_{R_{j}} & = & \textbf{v}, 
%\end{array}
%$$
%for some sequences $(\omega_{R_{j}})_{j}$ and $(\textbf{v}_{R_{j}})_{j}$. Thus, given that the map $A:\Omega\rightarrow\mathbb{C}^{2l}$ is continuous, for each $\epsilon>0$ and each $n\in\mathbb{Z}$, there exists $J\in\mathbb{N}$ such that $|n|\le R_J$ and for each $|r|\le R_J$, 
%$$
%\left\|A_{r}(\omega)\textbf{v}\right\| < 1+\epsilon.
%$$
\end{proof}

\begin{4}
\label{prop.hip.aberto}
Let $z\in\C$ and let $(T, A(z,\cdot))$ be the Jacobi cocycle defined over the compact metric space $\Omega$ by the law~\eqref{eq.cociclo.din}, with $T: \Omega \rightarrow \Omega$ a homeomorphism and $D,V:\Omega\rightarrow M(l,\mathbb{R})$ continuous maps. %, with $D(\omega)$ invertible for each $\omega\in\Omega$.
If $(T, A(z,\cdot)) \in \mathcal{UG}$, then there exists $\eta>0$ such that for each $z^\prime\in B(z;\eta)$, $(T, A(z^\prime,\cdot)) \in \mathcal{UG}$. Consequently, the condition $(T, A(z,\cdot))\in\mathcal{UG}$ is open with respect to $z\in\C$.
\end{4}
\begin{proof}
  Let $z\in\C$ be such that $(T, A(z,\cdot)) \in \mathcal{UG}$. Then, it follows from Definition~\ref{def.cresc.uni} that there exists $n_0=n_0(z)\in\N$ such that for each $(\omega,\textbf{v})\in\Omega\times\mathbb{S}^{2l-1}$,
  \[\Vert A_{n}(z,\omega)\textbf{v}\Vert>3/2,\;\forall\; n\ge n_0,\qquad\textrm{otherwise}\qquad \Vert A_{-n}(z,\omega)\textbf{v}\Vert>3/2,\;\forall\; n\ge n_0,
  \]
  and so, for each $n\ge n_0$ and each $(\omega,\textbf{v})\in\Omega\times\mathbb{S}^{2l-1}$,
  \[B_n(z,\omega,\textbf{v}):=\max\{\Vert A_{n}(z,\omega)\textbf{v}\Vert,\Vert A_{-n}(z,\omega)\textbf{v}\Vert\}>3/2.
  \]
  Fix $n\ge n_0$. %Assume, without loss of generality, that there exists $n\ge n_0$ such that $\Vert A_{n}(z,\omega)\textbf{v}\Vert>1,5$; fix such $n$.
  Then, by Definition~\eqref{def.mat.trans.erg} and by the fact that %$\Omega$ is a compact metric space and that
  $D, V : \Omega\rightarrow M (l, \R)$ are uniformly continuous maps (with $D(\omega)$ invertible for each $\omega\in \Omega$), it follows that there exists $\eta=\eta(z,n)>0$ such that for each $z^\prime\in B(z;\eta)$ and each $(\omega,\textbf{v})\in\Omega\times\mathbb{S}^{2l-1}$,
  %\[\Vert A_{n}(z^\prime,\omega)\textbf{v}\Vert>1,5;\]
\[B_n(z^\prime,\omega,\textbf{v})>3/2;\]
namely, $\{z\in\C\mid\sup_{(\omega,\textbf{v})\in\Omega\times\mathbb{S}^{2l-1}}B_n(z,\omega,\textbf{v})>3/2\}$ is open, given that
\[\C\ni z\mapsto\sup_{(\omega,\textbf{v})\in\Omega\times\mathbb{S}^{2l-1}}B_n(z,\omega,\textbf{v})\in\R\] is a lower-semicontinuous function.

Let $z^\prime\in B(z;\eta)$.  The result follows now for Lemma~\ref{lema.carac.omega.v}; let $\epsilon=1/2$, $R=n$ and $r\in\{-n,n\}$ in the statement of the lemma in order to conclude that $(T,A(z^\prime,\cdot))\in\mathcal{UG}$.
 \end{proof}

\section{Exponential Dichotomy}
\label{expdic}
\zerarcounters

In this section, we present the proofs of Theorems~\ref{prop.espec.2} and~\ref{prop.espec.3}. %The idea now is to use Proposition~\ref{prop.caracter} to prove the exponential dichotomy for the dynamically defined operator $H_{\omega}$.

\subsection{Proof of Theorem~\ref{prop.espec.2}}%the set inclusion $\rho(H_\omega)\subset\{z\in\C\mid (T,A_z)\in\mathcal{UG}\}$}

Firstly, we recall the Weyl Criterion, used for characterizing the spectrum of a linear operator defined in a Hilbert space (see, for example,~\cite{cesar09} for a proof).

\begin{4}[Weyl Criterion]
\label{crit.weyl}
Let $H: \dom(H) \subset \mathcal{H} \rightarrow \mathcal{H}$ be a linear operator defined in a Hilbert space $\mathcal{H}$ and let $z \in \mathbb{C}$. If there exists a sequence $(\textbf{v}_{n})_{n \in \mathbb{N}}$ of unitary vectors in $\dom(H)$ such that 
\[\lim_{n \rightarrow \infty} (H - z)\textbf{v}_{n} = \textbf{0},
\]
then $z \in \sigma(H)$.
\end{4}

\begin{3}
\label{lema.weyl}
Let $c_{00}(\mathbb{Z}, \mathbb{C}^{l}):=\{\textbf{u}\in(\mathbb{C}^{l})^{\mathbb{Z}}\mid \exists N\in\mathbb{Z}_0^+$ such that $\forall\, |j|>N, \textbf{u}_j=0\}$, and let $(H_{\omega})_{\omega}$ be the family of self-adjoint operators given by \eqref{eq.ope.din}, with $\dom(H_{\omega}) = l^{2}(\mathbb{Z}, \mathbb{C}^{l})$ for every $\omega \in \Omega$. Then, for each $\omega \in \Omega$, $x \in \sigma(H_{\omega})$ if, and only if for every $\epsilon > 0$, there exists a unitary vector $\textbf{u} \in c_{00}(\mathbb{Z}, \mathbb{C}^{l})$ such that
$$
\left\| (H_{\omega} - x)\textbf{u} \right\| < \epsilon.
$$
\end{3}
\begin{proof}
Let $\omega\in\Omega$, $\epsilon > 0$, and set $\delta > 0$ such that $0<\frac{2A\delta}{4A - \delta} \leq \epsilon$, where $A:=\sup\{\Vert H_\omega-x\Vert\mid x\in\sigma(H_\omega)\}<\infty$ (since $H_\omega$ is a bounded operator, $\sigma(H_\omega)\subset\R$ is bounded).  Let also $x \in \sigma(H_{\omega})$; it follows from Proposition \ref{crit.weyl} that there exists a unitary vector $\textbf{v} \in l^{2}(\mathbb{Z}, \mathbb{C}^{l})$ such that
$$
\left\| (H_{\omega} - x)    \textbf{v} \right\| < \frac{\delta}{4}.
$$  
For each $q \in \mathbb{N}$, define the vector  $\textbf{v}^{(q)} \in l^{2}(\mathbb{Z}, \mathbb{C}^{l})$ by %(\mathbb{C}^{l})^{\mathbb{Z}}$ by
$$
\textbf{v}^{(q)}_{n} := 
\left\{
\begin{array}{lll}
\textbf{v}_{n}, & n \in [-q, q], \\
&  \\
\textbf{0}, & n \notin [-q, q]. \\
\end{array}\right.
$$

Now, %since $H_\omega-x\mathbb{I}$ is a bounded operator,
there exists  $L \in \mathbb{N}$ such that %for each $q>L$,
$\left\| \textbf{v}^{(L)} - \textbf{v}\right\| < \frac{\delta}{4A}$ and so, 
$$
\left\| (H_{\omega} - x)    \textbf{v}^{(L)} \right\| < \frac{\delta}{2}.
$$
%(given that $H_\omega-x\mathbb{I}$ is a bounded operator).

The result follows now by taking $\textbf{u} := \frac{\textbf{v}^{(L)}}{\left\|\textbf{v}^{(L)}\right\|}$. Namely, %by the triangular inequlity,
$$
1 = \left\| \textbf{v} \right\| \leq  \left\| \textbf{v}^{(L)} \right\| +  \left\| \textbf{v}^{(L)} - \textbf{v}\right\| \leq  \left\| \textbf{v}^{(L)} \right\| + \frac{\delta}{4A},
$$
from which follows that $\left\| \textbf{v}^{(L)} \right\| \geq \frac{4A - \delta}{4A}$; then, 
$$
\left\| (H_{\omega} - x) \textbf{u} \right\| = \frac{\left\| (H_{\omega} - x) \textbf{v}^{(L)} \right\|}{\left\| \textbf{v}^{(L)} \right\|} \leq \left\| (H_{\omega} - x) \textbf{v}^{(L)} \right\| \left(\frac{4A}{4A - \delta}\right) < \frac{\delta}{2} \left(\frac{4A}{4A - \delta}\right) \leq \epsilon.
$$

The converse follows from Proposition~\ref{crit.weyl}.
\end{proof}

\begin{7}[Constancy of the spectrum]\label{constancy} Let $(H_{\omega})_{\omega}$ be the family of bounded dynamically defined operators given by~\eqref{eq.ope.din}, where $\Omega$ is a compact metric space, $T: \Omega \rightarrow \Omega$ is a minimal homeomorphism and $D,V:\Omega\rightarrow M(l,\mathbb{R})$ are continuous maps, with $D(\omega)$ invertible for each $\omega\in\Omega$. Then, for each $\omega_1,\omega_2 \in \Omega$, 
\[\sigma(H_{\omega_1})=\sigma(H_{\omega_2}).\]
\end{7}
\begin{proof}
Fix arbitrary~$\omega_1,\omega_2\in\Omega$.  Since~$(T,\Omega)$ is minimal, there exists a sequence~$(a_s)_{s\ge 1}$ such that~$\lim_{s\to\infty}T^{a_s}\omega_2=\omega_1$. By the continuity of~$D$ and $V$ at~$\omega_1$,~$H_{T^{a_s}\omega_2}$ converges strongly to~$H_{\omega_1}$ as~$s\to\infty$. Thus, by Corollary 10.2.2 in~\cite{cesar09} and the fact that, for each~$j\in\mathbb{Z}$ and each~$\eta\in\Omega$,~$H_{T^j\eta}$ and~$H_{\eta}$ have the same spectrum (since they are unitarily equivalent), one has
\[\sigma(H_{\omega_1})\subset\overline{\bigcup_{s\ge 1}\sigma(H_{T^{a_s}\omega_2})}=\sigma(H_{\omega_2}).\]
By interchanging the roles of $\omega_1$ and $\omega_2$ in the previous argument, one concludes that~$\sigma(H_{\omega_2})\subset\sigma(H_{\omega_1})$.
\end{proof}

\begin{remark}\label{rconstancy} We note that if the dynamical system $(T,\Omega)$ is topological transitive (that is, if there exists $\omega_0\in\Omega$ such that $\overline{\{T^j\omega_0\mid j\in\mathbb{Z}\}}=\Omega$), then it follows from the proof of Theorem~\ref{constancy} that for each $\omega\in\Omega$,
\[\sigma(H_{\omega})\subset\sigma(H_{\omega_0})\](compare this with Proposition~1 in~\cite{alkorn}).    
  \end{remark}

\begin{proof4}
 % For each $\omega\in\Omega$, let $V = (V_n)_{n\in\mathbb{Z}}$ and $D = (D_n)_{n\in\mathbb{Z}}$ be the bilateral sequences of real and symmetric $l\times l$ matrices that define $H_\omega$. Denote by $S:l_\infty(\mathbb{Z}; M (l, \mathbb{C}))\rightarrow l_\infty(\mathbb{Z}; M (l, \mathbb{C}))$ the shift operator and let $\Lambda:=\Lambda_1\times \Lambda_2$, with $\Lambda_1:=\overline{\{S^n(D)\mid n\in\mathbb{N}\}}$ and $\Lambda_2:=\overline{\{S^n(V)\mid n\in\mathbb{N}\}}$ (the closures are taken with respect to the topology of the uniform convergence over $l_\infty(\mathbb{Z}; M (l, \mathbb{C}))$).
%
 % Let also $g:\Lambda_{i}:\rightarrow SL(2l,\mathbb{C})$, $i=1,2$, be given by the law $g(B)=B_0$ (that is, $g$ is the evaluation map of the bilateral sequence $B\in\Lambda_{i}$ at entry $n=0$). Then, for each $z\in\mathbb{C}$, one defines the cocycle
%\begin{equation*}
%\label{eq.cociclo.din}
%\begin{array}{lllll}
%(S,P(z,\cdot)): & \Lambda \times \mathbb{C}^{2l} & \rightarrow & \Lambda \times \mathbb{C}^{2l}\\ 
%&(\lambda, \textbf{v}) & \mapsto & (S(\lambda), P(z,\lambda)\textbf{v}),
%\end{array}
%\end{equation*}  
%where the map $P(z,\cdot): \Omega \rightarrow SL(2l, \mathbb{C})$ is given by the law
%\begin{equation*}
%\label{def.func.mat.coci}
%P(z,\lambda)
%:=
%\left[
%\begin{array}{cc}
%g(B)^{-1}(z - g(C))    & - g(B)^{-1} \\
%& \\
%g(B) & 0
%\end{array}
%\right],
%\end{equation*}
%with $\lambda = (B, C)\in\Lambda_1\times\Lambda_2$.

  %Naturally, if the cocycle defined by~\eqref{eq.cociclo.din} is such that $(T, A(z,\cdot))\notin \mathcal{UG}$,
  Let us assume that  $(T,A(z,\cdot))\notin \mathcal{UG}$; then, it follows from Proposition~\ref{prop.caracter} that there exist $\omega \in \Omega$ and $\textbf{v} \in \mathbb{S}^{2l-1}$  such that for each $n \in \mathbb{Z}$,
\begin{equation}
  \label{sequlim}
\left\|A_{n}(z,\omega)\textbf{v}\right\| \leq 1,
\end{equation}
where $A_n(z,\omega)$ is given by~\eqref{def.mat.trans.erg}. %, with $P(z,\lambda)$ replacing $A(z,\omega)$.

Let $(\textbf{u}_{n})_{n \in \mathbb{Z}} \in (\mathbb{C}^{l})^{\mathbb{Z}}$ be the bilateral sequence given by the law 
$$
\left[
\begin{array}{c}
\textbf{u}_{n+1}\\
D(T^n\omega)\textbf{u}_{n}
\end{array}\right]
=
A_{n}(z,\omega)\textbf{v},
$$
with $u_1:=v_1$ and $u_0:=(D(\omega))^{-1}v_0$ (here, $v_1,v_0\in\mathbb{C}^l$ are the components of $\textbf{v}$).

%Denote by $H_\lambda$ the matrix-valued Jacobi operator~\eqref{eq.ope.din.jacobi} associated with $\lambda=(B,C)$.
Then, $\textbf{u}=(\textbf{u}_{n})_n$ is a solution to the eigenvalue equation of $H_{\omega}$ at $z$ and by~\eqref{sequlim}, $\left\|\textbf{u}\right\|_{\infty}<\infty$. %(namely, given that $D:\Omega\rightarrow M(l,\mathbb{R})$ is continuous and $\Omega$ is compact, it follows that $M:\Vert D\Vert_\infty<.
If $(\textbf{u}_{n}) \in l^{2}(\mathbb{Z}, \mathbb{C}^{l})$, then $z \in \sigma(H_{\omega})$ and we are done. Otherwise, for each $L \in \mathbb{N}$, one may define $\textbf{u}^{(L)}$ by the law
\begin{eqnarray}\label{seqWeyl}
\textbf{u}^{(L)}_{n} = 
\left\{
\begin{array}{cc}
\textbf{u}_{n}, &  \left| n \right| \leq L,  \\
&\\
\textbf{0}, & \left|n\right| > L.
\end{array}\right.
\end{eqnarray}

If $n\notin\{-L-1,-L,-L+1,L-1,L,L+1\}$, since $\textbf{u}$ is a solution to the eigenvalue equation of $H_\omega$ at $z$, one has 
$$
[(H_{\omega} - z\mathbb{I})\textbf{u}^{(L)}]_{n} = \textbf{0}.
$$ 
Specifically at the entries $n = \pm (L - 1), \pm L, \pm (L + 1)$, one has 
$$
\left\{
\begin{array}{lll}
  \left[(H_{\omega} - z\mathbb{I})\textbf{u}^{(L)}\right]_{\pm(L - 1)} & = & D(T^{\pm(L - 2)}\omega)\textbf{u}_{\pm(L - 2)} + D(T^{\pm(L - 1)}\omega)\textbf{u}_{\pm L} \\
  &&\\
  &+& (z\mathbb{I} - V(T^{\pm(L - 1)}\omega))\textbf{u}_{\pm(L - 1)}, \\
& & \\
\left[(H_{\omega} - z\mathbb{I})\textbf{u}^{(L)}\right]_{\pm L} & = & D(T^{\pm(L - 1)}\omega) \textbf{u}_{\pm(L - 1)} +  D(T^{\pm L}\omega)\textbf{u}_{\pm L}, \\
& & \\
\left[(H_{\omega} - z\mathbb{I})\textbf{u}^{(L)}\right]_{\pm(L + 1)} & = &  D(T^{\pm L}\omega) \textbf{u}_{\pm L}.
\end{array}\right. .
$$
%it is possible to obtain analogous results for $n = - L - 1, -L, - L + 1$.
Thus, %, and since for each $|n|>L+1$,
%$$
%[(H_{\lambda} - z\mathbb{I})\textbf{u}^{(L)}]_{n} = \textbf{0},
%$$
one concludes that there exists a constant $K$ such that for each $L \in \mathbb{N}$, 
$$
\left\|(H_{\omega} - z\mathbb{I})\textbf{u}^{(L)}\right\| \leq K
$$
(note that $(V_{n})_{n \in \mathbb{Z}}$ and $(D_{n})_{n \in \mathbb{Z}}$ are uniformly bounded bilateral sequences, and that for each $L\in\N$, $\Vert\textbf{u}^{(L)}\Vert_\infty\le\Vert\textbf{u}\Vert_\infty<\infty$).

Now, by setting $\textbf{w}^{(L)} := \frac{\textbf{u}^{(L)}}{\left\| \textbf{u}^{(L)}  \right\|}$, it follows from the last inequality that %by $\left\| \textbf{u}^{(L)} \right\|$, we have that 
$$
\left\|(H_{\omega} - z\mathbb{I})\textbf{w}^{(L)}\right\| \leq \frac{K}{\left\| \textbf{u}^{(L)}  \right\|}.
$$ 

Since, by hypothesis, $\lim_{L \rightarrow \infty}\left\| \textbf{u}^{(L)}  \right\| = \infty$ (given that $\textbf{u}^{(L)}\notin l^2(\mathbb{Z},\mathbb{C}^l)$), it follows from Lemma~\ref{lema.weyl} that $z \in \sigma(H_{\omega})$.

Therefore, if $(T,A(z,\cdot))\notin\mathcal{UG}$, then there exists $\omega\in\Omega$ such that $z \in \sigma(H_{\omega})$. The result is now a direct consequence of Theorem~\ref{constancy}.
\end{proof4}

\begin{remark}\label{Rconstancy} It follows from the proof of Theorem~\ref{prop.espec.2} that if $T$ is assumed to be a homeomorphism (without any dynamical hypothesis over $(T,\Omega)$, such as minimality or transitivity), then
  \[\bigcap_{\omega\in\Omega}\rho(H_{\omega})\subseteq\{z\in\mathbb{C}\mid(T,A(z,\cdot))\in\mathcal{UG}\}.\]
  Moreover, if $(T,\Omega)$ is assumed to be transitive, then it follows from Remark~\ref{rconstancy} and from the argument in the previous paragraph that
  \[\rho(H_{\omega_0})\subseteq\{z\in\mathbb{C}\mid(T,A(z,\cdot))\in\mathcal{UG}\},\]
where $\omega_0\in\Omega$ is such that  $\overline{\{T^j\omega_0\mid j\in\mathbb{Z}\}}=\Omega$.
\end{remark}
%%%%%%%%%%%%%%%%%%%%%%%%%%%%%%%%%%%%%%%%%%%%%%%%%%%%%%%%%%%%%%%%%%%%%%%%%%%%%%%%%%%%%%%%%%%%%%%%%%%%%%%%%%%%%%%%%%%%%%%%%%%%%%%%%%%%%%%%%%%%%%%%%%%%%%%%%%%%%%%%%%%%%%%%%%%%%%%%%%%%%%%%%%%%%%%%%%%%%%%%%%%%%%%%%%%%%%%%%%%%%%%%%%%%%%%%%%%%%%%%%%%%%%%%%%%%%%

\subsection{Proof of Theorem~\ref{prop.espec.3}}

Now we prove the converse of Theorem~\ref{prop.espec.2}. %, concluding the proof of Theorem~\ref{maintheo}.
The idea is to show that if there exist $2l$ linearly independent solutions to the eigenvalue equation of $H_\omega$ at $z$ %, say $\{\textfb{u}_1,\ldots,\textfb{u}_{2l}\}$,
such that half of them belong to $l^2(\mathbb{N},\mathbb{C}^{l})$ and the other half belong to $l^2(\mathbb{Z}_-,\mathbb{C}^{l})$, then $z\in\rho(H_\omega)$. In order to prove this result, % belongs to the resolvent set. For this,
one needs to write explicitly the resolvent operator on its integral form, the so-called Green Function (Proposition \ref{porp.def.green}). Some preparation is required.

%In what follows, we fix $\omega\in\Omega$ and denote $H_\omega$, $D(T^n\omega)$ and $V(T^n\omega)$ ($n\in\mathbb{Z}$), respectively, by $H$, $D_n$ and $V_n$. 

We begin with Green Formula. Let $\textbf{u}, \textbf{v} \in (\mathbb{C}^{l})^{\mathbb{Z}}$; then, for every integers $n > m$, 
\begin{equation}
\label{eq.green.1}
\sum^{n}_{k = m} \left\langle (H\textbf{u})_{k}, \bar{\textbf{v}}_{k} \right\rangle_{\mathbb{C}^{l}} - \left\langle (H\textbf{v})_{k}, \bar{\textbf{u}}_{k} \right\rangle_{\mathbb{C}^{l}} = W_{[\textbf{u}, \textbf{v}]}(n + 1) - W_{[\textbf{u}, \textbf{v}]}(m),
\end{equation}
where $W$ is the Wronskian of $\textbf{u}$ and $\textbf{v}$, given by
\begin{equation}
\label{eq.wrons1}
W_{[\textbf{u}, \textbf{v}]}(n) :=  \left\langle D_{n - 1}\textbf{u}_{n}, \bar{\textbf{v}}_{n - 1} \right\rangle_{\mathbb{C}^{l}} - \left\langle D_{n - 1} \textbf{v}_{n}, \bar{\textbf{u}}_{n - 1}  \right\rangle_{\mathbb{C}^{l}}.
\end{equation} 

If $(A_{n})$ and $(B_{n})$ are sequences in $M(l,\mathbb{C})$, one obtain, by
applying the operator to each of their columns, the following version of Green Formula for
matrices: 
\begin{equation}
\label{wronski.matriz}
\sum^{n}_{k = m} A^{\ast}_{k} H(B)_{k} - H(A)_{k}^{\ast} B_{k}  =  W_{[A, B]}(n + 1) - W_{[A, B]}(m),
\end{equation}
with
$$
W_{[A, B]}(m) = (A^{\ast}_{m - 1}D_{m - 1}B_{m} - A^{\ast}_{m}D_{m - 1}B_{m - 1}).
$$

In the specific case that $\textbf{u}_n$ and $\textbf{v}_n$ are solutions to the same eigenvalue equation, one obtains the constancy of the Wronskian.

\begin{3}[Constancy of the Wronskian]
\label{lema.const.wronsk}
Let $H$ be the operator given by \eqref{eq.ope.din.jacobi}. If $\textbf{u}, \textbf{v} \in (\mathbb{C}^{l})^{\mathbb{Z}}$ satisfy the eigenvalue equation for $H$ %\textbf{u} = z \textbf{u}$ and $H \textbf{v} = z \textbf{v}$ for some
at $z \in \mathbb{C}$, then for each $m,n\in\mathbb{Z}$,
\[W_{[\textbf{u}, \textbf{v}]}(m)=W_{[\textbf{u}, \textbf{v}]}(n).
\]
\end{3}
\begin{proof}
It is enough to apply the Green formula for $\textbf{u}$ and $\textbf{v}$. Namely,
$$
\sum^{n}_{k = m} \left\langle (H\textbf{u})_{k}, \bar{\textbf{v}} \right\rangle_{\mathbb{C}^{l}} - \left\langle (H\textbf{v})_{k},\bar{\textbf{u}}_{k} \right\rangle_{\mathbb{C}^{l}}  =  \sum^{n}_{k = m} \left\langle z \textbf{u}_{k},\bar{\textbf{v}}_{k} \right\rangle_{\mathbb{C}^{l}} - \left\langle z \textbf{v}_{k}, \bar{\textbf{u}}_{k} \right\rangle_{\mathbb{C}^{l}} = 0,
$$
from which follows that
$$
W_{[\textbf{u}, \textbf{v}]}(n + 1) - W_{[\textbf{u}, \textbf{v}]}(m) = 0. 
$$
\end{proof}
 
We want to define the resolvent operator in terms of the solutions to the eigenvalue equation. 

\begin{3}
\label{lema.green.nd}
Let $H$  be given by \eqref{eq.ope.din.jacobi} and let $(F_{n}^{(+)})_{n\in\mathbb{Z}}, (F_{n}^{(-)})_{n\in\mathbb{Z}}$ be bilateral sequences in $M(l,\mathbb{C})$ whose columns are linearly independent solutions to the eigenvalue equation of $H$ at $z \in \mathbb{C}$. Suppose that 
$$
\lim_{n \rightarrow \pm \infty} \left\| F^{(\pm)}_{n} \right\| = 0
$$
and set
$$
Q:= W_{[F^{(+)}, F^{(-)}]}(0) = (F^{(+)}_{0})^{\ast}D_{0}F^{(-)}_{1} - (F^{(+)}_{1})^{\ast}D_{0}F^{(-)}_{0}.
$$ 
Then, for every $n \in \mathbb{Z}$, 
\begin{eqnarray*}
(a) & F^{(+)}_{n}Q^{-1}(F^{(-)}_{n})^{\ast} - F^{(-)}_{n}(Q^{\ast})^{-1}(F^{(+)}_{n})^{\ast}=0;\\
(b) & F^{(+)}_{n}Q^{-1}(F^{(-)}_{n + 1})^{\ast} - F^{(-)}_{n}(Q^{\ast})^{-1}(F^{(+)}_{n + 1})^{\ast}  =  D^{-1}_{n};\\ 
(c) & F^{(+)}_{n + 1}Q^{-1}(F^{(-)}_{n})^{\ast} - F^{(-)}_{n + 1}(Q^{\ast})^{-1}(F^{(+)}_{n})^{\ast} = - D^{-1}_{n}. 
\end{eqnarray*}
\end{3}

\begin{proof}
By letting $A_{n}=B_{n} = F^{(\pm)}_{n}$ in Green Formula \eqref{wronski.matriz}, it follows that there exists a constant $C_\pm\in M(l,\mathbb{C})$ such that for every $n \in \mathbb{Z}$,
\[
W_{[F^{(\pm)}, F^{(\pm)}]}(n) = ((F^{(\pm)}_{n - 1})^{\ast}D_{n - 1}F^{(\pm)}_{n} - (F_{n}^{(\pm)})^{\ast}D_{n - 1}F_{n - 1}^{(\pm)}) = C_\pm;
\]
since $\lim_{n \rightarrow \pm \infty} \left\| F^{(\pm)}_{n} \right\| = 0$, one has
\[
\lim_{n \rightarrow \infty} W_{[F^{(\pm)}, F^{(\pm)}]}(n) = 0,
\]
and then, $C_\pm = 0$.

Now, if one applies equation~\eqref{wronski.matriz} to the pairs $(F^{(+)},F^{(-)}), (F^{(-)},F^{(+)})$ and $(F^{(-)},F^{(-)})$, one obtains from the constancy of the Wronskian, for each $n \in \mathbb{Z}$, the system
$$
\left\{
\begin{array}{lll}
(F^{(+)}_{n})^{\ast}D_{n}F^{(+)}_{n + 1} - (F^{(+)}_{n + 1})^{\ast}D_{n}F^{(+)}_{n} & = & 0, \\ 
&&\\
(F^{(+)}_{n})^{\ast}D_{n}F^{(-)}_{n + 1} - (F^{(+)}_{n + 1})^{\ast}D_{n}F^{(-)}_{n} & = & Q, \\
&&\\
(F^{(-)}_{n})^{\ast}D_{n}F^{(+)}_{n + 1} - (F^{(-)}_{n + 1})^{\ast}D_{n}F^{(+)}_{n} & = & -Q^{\ast}, \\
&&\\
(F^{(-)}_{n})^{\ast}D_{n}F^{(-)}_{n + 1} - (F^{(-)}_{n + 1})^{\ast}D_{n}F^{(-)}_{n} & = & 0,
\end{array}\right.
$$
which can be written in the form 
\begin{equation}
\label{equa.matri.simpl}
\left[
\begin{array}{cc}
(F^{(+)}_{n})^{\ast}_{n} & (F^{(+)}_{n})^{\ast}_{n + 1}\\
(F^{(-)}_{n})^{\ast}_{n} & (F^{(-)}_{n})^{\ast}_{n + 1}
\end{array}
\right]
\left[
\begin{array}{cc}
D_{n} & 0 \\
0 & D_{n}
\end{array}
\right]
\mathbb{J}
\left[
\begin{array}{cc}
F^{(+)}_{n} & F^{(-)}_{n} \\
F^{(+)}_{n + 1} & F^{(-)}_{n + 1}
\end{array}
\right]
=
\mathbb{J}
\left[
\begin{array}{cc}
Q & 0 \\
0 & Q^{\ast}
\end{array}
\right],
\end{equation}
with
$$
\mathbb{J}
=
\left[
\begin{array}{cc}
0 & \mathbb{I} \\
- \mathbb{I} & 0
\end{array}
\right].
$$

Since $\mathbb{J}^{-1} = -\mathbb{J}$ and 
$$
\left[
\begin{array}{cc}
Q & 0 \\
0 & Q^{\ast}
\end{array}
\right]^{-1} = \left[
\begin{array}{cc}
Q^{-1} & 0 \\
0 & (Q^{\ast})^{-1}
\end{array}
\right],
$$
by multiplying to the left both members of identity \eqref{equa.matri.simpl} by  $\left[
\begin{array}{cc}
Q^{-1} & 0 \\
0 & (Q^{\ast})^{-1}
\end{array}
\right]\mathbb{J}^{-1}$, it follows that
$$
\left(
\left[
\begin{array}{cc}
F^{(+)}_{n} & F^{(-)}_{n} \\
F^{(+)}_{n + 1} & F^{(-)}_{n + 1}
\end{array}
\right]
\right)
\left(
\left[
\begin{array}{cc}
Q^{-1} & 0 \\
0 & (Q^{\ast})^{-1}
\end{array}
\right]
\mathbb{J}^{-1}
\left[
\begin{array}{cc}
(F^{(+)}_{n})^{\ast} & (F^{(+)}_{n + 1})^{\ast}\\
(F^{(-)}_{n})^{\ast} & (F^{(-)}_{n + 1})^{\ast}
\end{array}
\right]
\left[
\begin{array}{cc}
D_{n} & 0 \\
0 & D_{n}
\end{array}
\right]
\mathbb{J}
\right)
=
\left[
\begin{array}{cc}
\mathbb{I} & 0 \\
0 & \mathbb{I} 
\end{array}
\right],
$$
since one-sided finite matrix inverses are two-sided; in other words,
$$
\left[
\begin{array}{cc}
Q^{-1} & 0 \\
0 & (Q^{\ast})^{-1}
\end{array}
\right]
\mathbb{J}^{-1}
\left[
\begin{array}{cc}
(F^{(+)}_{n})^{\ast} & (F^{(+)}_{n + 1})^{\ast}\\
(F^{(-)}_{n})^{\ast} & (F^{(-)}_{n + 1})^{\ast}
\end{array}
\right]
\left[
\begin{array}{cc}
D_{n} & 0 \\
0 & D_{n}
\end{array}
\right]
\mathbb{J}
\;\;\;\; \mbox{and} \;\;\;\;
\left[
\begin{array}{cc}
F^{(+)}_{n} & F^{(-)}_{n} \\
F^{(+)}_{n + 1} & F^{(-)}_{n + 1}
\end{array}
\right]
$$
commute. Such identity can be written in the form
$$
\left\{
\begin{array}{lll}
F^{(+)}_{n}Q^{-1}(F^{(-)}_{n})^{\ast} - F^{(-)}_{n}(Q^{\ast})^{-1}(F^{(+)}_{n})^{\ast} & = & 0, \\ 
&&\\
F^{(+)}_{n}Q^{-1}(F^{(-)}_{n + 1})^{\ast} - F^{(-)}_{n}(Q^{\ast})^{-1}(F^{(+)}_{n + 1})^{\ast} & = & D^{-1}_{n}, \\
&&\\
F^{(+)}_{n + 1}Q^{-1}(F^{(-)}_{n})^{\ast} - F^{(-)}_{n + 1}(Q^{\ast})^{-1}(F^{(+)}_{n})^{\ast} & = & - D^{-1}_{n}, \\
&&\\
F^{(+)}_{n + 1}Q^{-1}(F^{(-)}_{n + 1})^{\ast} - F^{(-)}_{n + 1}(Q^{\ast})^{-1}(F^{(+)}_{n + 1})^{\ast} & = & 0.
\end{array}\right.
$$
\end{proof}

One can obtain from these relations the Green Function of $H$.%Now we can define the reolvent operator by the Green formula: 

\begin{4}
\label{porp.def.green}
Let $H$ be the operator given by \eqref{eq.ope.din.jacobi}, let $z\in\mathbb{C}$  and let $F_{n}^{(+)}, F_{n}^{(-)}$ be sequences of matrices in $M(l,\mathbb{C})$ such that $\lim_{n\to\pm\infty}\Vert F_n^{(\pm)}\Vert=0$ and whose columns are linearly independent solutions to the eigenvalue equation of $H$ at $z$. %Let $A:=\{z\in\mathbb{C}\mid$ there exist constants $C_1,C_2>0$ and $\lambda>1$ such that for each $n\in\Z$, $\Vert F_n^+(z)\Vert\le C_1\lambda^{-n}$, $\Vert F_n^-(z)\Vert\le C_2\lambda^{n}\}$.  
%$$
%\sum^{\pm \infty}_{n = 0} \left\|F_{n}^{(\pm)}\right\|^{2} < \infty.
%$$
 Set, for each $p, q \in \mathbb{Z}$, % and each $z \in \mathbb{C}\setminus\mathbb{R}$,%For $p, q \in \mathbb{Z}$, define
\begin{equation}\label{eq.resol.green}
G(p, q; z) = 
\left\{
\begin{array}{ll}
- F_{p}^{(-)}(z) (Q^{\ast})^{-1} (F_{q}^{(+)})^{\ast}(z), & p \leq q, \\
&\\
- F^{(+)}_{p}(z) Q^{-1} (F^{(-)}_{q})^{\ast}(z), & p > q.
\end{array}\right.
\end{equation}
Then, for each $\textbf{u} \in (\mathbb{C}^{l})^{\mathbb{Z}}$,
\begin{equation*}
%\label{eq.resol.green}
\sum_{q} G(p, q; z)\textbf{u}_{q} = ((H - z)^{-1}\textbf{u})_{p};
\end{equation*}
that is, $G$ is the Green Function of the operator $H$ at $z$. %Moreover, $(H-z)^{-1}$ is bounded for each $z\in A$.

\end{4}
\begin{proof}
 % Let $\textbf{u}\in l^2(\mathbb{Z};\C^l)$; since %, for each $n\in\N$,
 % $F^{\pm}(z)\in l^{2}(\mathbb{Z}_\pm; \mathbb{C}^{l})$, it follows that  %$\textbf{u} \in l^{2}(\mathbb{Z}_{+}, \mathbb{C}^{l})$,
%\begin{eqnarray*}
%\left(\sum_{q} G(p, q; z)\textbf{u}_{q}\right)_{p \in \Z}\in l^{2}(\Z; \mathbb{C}^{l}).
%\end{eqnarray*}
%if $\textbf{u} \in l^{2}(\mathbb{Z}_{+};\mathbb{C}^{l})$.
One needs to prove that %We begin proving that
$$
(H - z)\left(\sum_{q} G(p, q; z)\textbf{u}_{q}\right)_{p \in \Z} = (\textbf{u}_{p})_{p \in \Z}.
$$

It follows from the definition of $G(p, q; z)$ that (we omit the dependence of $z$ in $G(p, q; z)$),
$$
\begin{array}{lll}
(H - z)\left(\sum_{q} G(p, q)\textbf{u}_{q}\right)_{p \in \Z}= &  & \left(\sum_{q} D_{p}G(p + 1, q)\textbf{u}_{q}\right)_{p \in \N}+\\
& & \\ 
& & \left(\sum_{q} D_{p - 1}G(p - 1, q)\textbf{u}_{q}\right)_{p \in \Z}+  \\
& & \\
& & \left(\sum_{q} (V_{p} - z)G(p, q)\textbf{u}_{q}\right)_{p \in \Z}.
\end{array}
$$

Let $p\in\Z$; one needs to consider the following cases.

\textbf{Case} $q < p - 1$:
\[
\begin{array}{lll}
& D_{p}G(p + 1, q)\textbf{u}_{q} + D_{p - 1}G(p - 1, q)\textbf{u}_{q} + (V_{p} - z)G(p, q)\textbf{u}_{q} = & \\
& & \\
& - D_{p}F_{p + 1}^{(+)}Q^{-1} (F^{(-)}_{q})^{\ast}\textbf{u}_{q} - D_{p - 1}F^{(+)}_{p - 1} D^{-1}_{0} (F^{(-)}_{q})^{\ast}\textbf{u}_{q} - (V_{p} - z) F_{p} Q^{-1}(F^{(-)}_{q})^{\ast}\textbf{u}_{q} = & \\
& & \\
 & - \left[ D_{p}F_{p + 1}^{(+)} + D_{p - 1}F_{p - 1}^{(+)} + (V_{p} - z) F_{p}^{(+)} \right]Q^{-1}(F^{(-)}_{q})^{\ast}\textbf{u}_{q}  = - \left[ 0 \right] Q^{-1}(F^{(-)}_{q})^{\ast}\textbf{u}_{q} = \textbf{0}. & 
\end{array}
\]

\textbf{Case} $q = p - 1$:
\[
\begin{array}{lll}
  &D_{p}G(p + 1, q)\textbf{u}_{q} + D_{p - 1}G(p - 1, q)\textbf{u}_{q} + (V_{p} - z)G(p, q)\textbf{u}_{q}=& \\
  & & \\
  &- D_{p}F^{(+)}_{p + 1} Q^{-1} (F^{(-)}_{q})^{\ast}\textbf{u}_{q} - D_{p - 1}F^{(-)}_{p - 1} (Q^{\ast})^{-1} (F^{(+)}_{q})^{\ast}\textbf{u}_{q} - (V_{p} - z) F^{(+)}_{p} Q^{-1}(F^{(-)}_{q})^{\ast}\textbf{u}_{q}=&\\
  & &  \\
&- D_{p}F^{(+)}_{p + 1} Q^{-1} (F^{(-)}_{q})^{\ast}\textbf{u}_{q} - D_{p - 1}F^{(+)}_{p - 1}Q^{-1}(F^{(-)}_{q})^{\ast}\textbf{u}_{q} - (V_{p} - z) F^{(+)}_{p}Q^{-1}(F^{(-)}_{q})^{\ast}\textbf{u}_{q}= &\\
&&\\
&- \left[ D_{p}F^{(+)}_{p + 1} + D_{p - 1}F^{(+)}_{p - 1} + (V_{p} - z) F^{(+)}_{p} \right] Q^{-1}(F^{(-)}_{q})^{\ast}\textbf{u}_{q} =- \left[ 0 \right] Q^{-1} \psi^{\ast}_{q} \textbf{u}_{q} = \textbf{0}&, 
\end{array}
\]
where it has been applied Lemma~\ref{lema.green.nd}-$(a)$ to the second identity.

\textbf{Case} $q=p$: by applying Lemma~\ref{lema.green.nd}-$(c)$, one has 
\[
D_{q}G(q + 1, q)\textbf{u}_{q} + D_{q - 1}G(q - 1, q)\textbf{u}_{q} + (V_{q} - z)G(q, q)\textbf{u}_{q} = \textbf{u}_{q}.
\]

\textbf{Case} $q \geq p + 1$:
\[
D_{p}G(p + 1, q)\textbf{u}_{q} + D_{p - 1}G(p - 1, q)\textbf{u}_{q} + (V_{p} - z)G(p, q)\textbf{u}_{q} = \textbf{0}.
\] 
%
%t remains to prove that for each $z\in A$ and each $\textbf{u} \in l^{2}(\mathbb{Z}, \mathbb{C}^{l})$,
%\[
%\left(\sum_{q} G(p, q; z)\textbf{u}_{q}\right)_p\in l^{2}(\mathbb{Z}, \mathbb{C}^{l}).\]
%
%We begin with the remark that for each $z\in A$, there exists a constant $D>0$ such that for each $p,q\in\Z$, $\Vert G(p, q; z)\Vert\le D\lambda^{-|p-q|}$; this is a direct consequence of the definition of $G(p,q;z)$. Now, one has, by Cauchy-Schwarz inequality and Fubini,
%\begin{eqnarray*}
 % \sum_p\left\Vert\sum_{q} G(p, q; z)\textbf{u}_{q}\right\Vert^2&\le&\sum_p\left(\sum_{q} \Vert G(p, q; z)\Vert\right)\left(\sum_q\Vert G(p,q;z)\Vert\Vert\textbf{u}_{q}\Vert^2\right)\\
 % &\le&D\,\left(\frac{1+\lambda}{1-\lambda}\right)\sum_q\Vert\textbf{u}_{q}\Vert^2\sum_p \Vert G(p,q;z)\Vert\le D^2\,\left(\frac{1+\lambda}{1-\lambda}\right)^2\sum_q\Vert\textbf{u}_{q}\Vert^2\\
 % &=&D^2\,\left(\frac{1+\lambda}{1-\lambda}\right)^2\Vert\textbf{u}\Vert^2,
%\end{eqnarray*}  
%and we are done.
\end{proof}

The next step consists in obtaining a minimal support for the trace of the spectral (matrix) measure of $H_\omega$, $\omega\in\Omega$.

We note that Green Functions $G_\omega(j, j; \cdot):\mathbb{C}_{+}\rightarrow M(l,\mathbb{C})$, $j\in\Z$, are matrix-valued Herglotz functions (that is, each $G_\omega(j, j; \cdot)$ is analytic and $\Im G_\omega(j,j;z)>0$ for each $z\in\mathbb{C}_+$; this is a consequence of the fact that $G_\omega$ is the integral kernel of $(H_\omega-z)^{-1}$ and that $\Im (H_\omega-z)^{-1}>0)$, from which follows that for $\kappa$-a.e.~$x\in\mathbb{R}$ (here, $\kappa$ stands for the Lebesgue measure on $\R$), 
\[
\lim_{y \downarrow 0}\Im G_\omega(j,j; x \pm iy)<\infty
\] 
(see~\cite{gesztesy97}). By Spectral Theorem one has, for each $\textbf{u} \in l^{2}(\mathbb{Z}; \mathbb{C}^{l})$, 
\begin{equation}
\label{eq.med.espec}
\left\langle  (H_\omega - z)^{-1} \textbf{u}, \textbf{u} \right\rangle = \int \frac{1}{x - z} d\mu^\omega_{\textbf{u}}(x),
\end{equation}
where $\mu_{\textbf{u}}(\cdot):=\langle\textbf{u},E^\omega(\cdot)\textbf{u}\rangle$ is a finite Borel measure and $E^\omega$ is the resolution of the identity of the operator $H_\omega$.

Let $H: \dom(H)\subset\mathcal{H} \rightarrow \mathcal{H}$ be a self-adjoint operator defined in a separable Hilbert space $\mathcal{H}$, and let $\mathcal{C} = \{\textbf{u}_{1}, \textbf{u}_{2}, \ldots, \textbf{u}_{k}\}\subset\mathcal{H}$. The cyclic subspace of $H$ spanned by $\mathcal{C}$ is the space
\begin{equation*}
%\label{eq.esp.cicl}
\mathcal{H}_{\mathcal{C}}:=\overline{\gerado\{\cup_{j=1}^k\{(H)^n(\textbf{u}_{j})\mid n\in\mathbb{N}\}\}}.%=l^{2}(\mathbb{N}; \mathbb{C}^{l}).%\overline{\gerado \left\{ \bigcup_{j = 1}^{k}\bigcup_{n \in \mathbb{N}} H^{n}(\textbf{u}_{j}) \right\}}.
\end{equation*}

One says that $\mathcal{C} = \{\textbf{u}_{1}, \textbf{u}_{2}, \ldots, \textbf{u}_{k}\}$ is a spectral basis for $H$ if the system $\mathcal{C}$ is linearly independent and $\mathcal{H}_{\mathcal{C}}=\mathcal{H}$.

In our setting, the $2l$ canonical vectors $(\textbf{e}_{\alpha,k})_{\alpha=0,1,k = 1,\ldots,l}$ in $(\mathbb{C}^{l})^{\Z}$, where $(\textbf{e}_{\alpha,k})_{n,j} = \delta_{\alpha,n}\delta_{j,k}$, $\alpha=0,1$, form a spectral basis for $H_\omega$ (see~\cite{carmona90}). The matrix
\begin{eqnarray*}
  \left(\begin{array}{cc}\mu^\omega_{0,0} & \mu^\omega_{0,1}\\
    \mu^\omega_{1,0}&\mu^\omega_{1,1}\end{array}\right),
\end{eqnarray*}
where $\mu^\omega_{\alpha,\beta}=(\mu^\omega_{\textbf{e}_{\alpha, i},\textbf{e}_{\beta, j}})_{1 \leq i,j \leq l}$, $\alpha,\beta=0,1$, is called the spectral (matrix) measure of $H_\omega$. % (we omit the dependence on $\omega$ from now on).
It can be shown (see~\cite{carmona90}) that the spectral type of the operator $H_\omega$ is given by $\mu^{tr}_\omega=\tr[\mu_{0,0}^\omega+\mu_{1,1}^\omega]$, in the sense that for each $\textbf{u}\in l^2(\Z;\C^l)$, $\mu^\omega_{\textbf{u}}$ is absolutely continuous with respect to $\mu^{tr}_\omega$.

%The next step consists in obtaining a characterization of a minimal support for the spectral measure $\mu^{tr}_\omega$.

\begin{1}[Minimal Support]
\label{def.sup.minimal}
One says that a set $S \subseteq \mathbb{R}$ is a minimal support for the positive and finite Borel measure $\mu$ if 
\begin{eqnarray*}\begin{array}{ll}
(i) & \mu (\mathbb{R} \setminus S) = 0; \\
(ii) & S_{0}\subset S,\;\;\; \mu(S_{0}) = 0\qquad \Longrightarrow\qquad \kappa(S_{0}) = 0.
\end{array}\end{eqnarray*} 
\end{1}

In other words, a minimal support for $\mu$ is a Borel set in which $\mu$ is concentrated and such that its subsets of zero measure necessarily have zero Lebesgue measure. Definition \ref{def.sup.minimal} induces an equivalence relation in $\mathcal{B}(\R)$ (the Borel $\sigma$-algebra of $\R$): %which the classes are given by the property
$$
S_{1} \sim S_{2} \Longleftrightarrow \kappa(S_{1} \Delta S_{2})  = \mu(S_{1} \Delta S_{2}) = 0,
$$
where $S_{1} \Delta S_{2}:= (S_{1} \setminus S_{2}) \cup (S_{2} \setminus S_{1}) $ is the symmetric difference of $S_{1}$ and $S_{2}$ (see Lemma 2.20 in \cite{gilbert84} for a proof of this statement).

\begin{4}[Minimal Supports for the Spectral Types of $H_\omega$]
\label{porp.sup.ac}
Let $H_\omega$, $\omega\in\Omega$, be given by \eqref{eq.ope.din} and let $G_\omega(z):=G_\omega(0,0;z)+G_\omega(1,1;z)$, where $G_\omega(j,j;z)$ is the Green Function of $H_\omega$ at $z\in\C$ evaluated at $j=0,1$, given by \eqref{eq.resol.green}. Then, %the set 
\begin{eqnarray*}
\Sigma_{ac}^\omega  &:=&  \{x \in \mathbb{R}\mid \exists \lim_{y \downarrow 0} G_\omega(x + iy),\; 0<\pi^{-1}\lim_{y \downarrow 0} \Im[\ln[G_\omega(x + iy)]]<\infty\},\\
%\vspace{2mm}
%is a minimal support for the absolutely continuous component of the spectral measure, %and %the set 
%$$
%\Sigma_{ac} := \bigcup^{l}_{j = 1} \Sigma^{\phi}_{ac, j}
%$$
%is a minimal support for the absolutely continuous component. Moreover, 
\Sigma_{s}^\omega &:=& \{x \in \mathbb{R}; \lim_{y \downarrow 0} \Im[\tr[G_\omega(x + iy)]] = \infty \}
\end{eqnarray*}
%and %is a minimal support for the singular component of the spectral measure, and
are minimal supports for the absolutely continuous and singular parts of the spectral measure $\mu^{tr}_\omega$, respectively. Moreover,
\[
\Sigma^\omega:=\Sigma_{ac}^\omega\cup\Sigma_{s}^\omega\]
is a minimal support for the spectral measure $\mu^{tr}_\omega$.
\end{4}
\begin{proof}
  See Theorem~6.1 in~\cite{gesztesy97} for details.
  \end{proof}

We are now ready to prove the converse of Theorem~\ref{prop.espec.2}.
%\begin{7}
%\label{prop.espec.3}
%Let $(H_{\omega})_{\omega}$ be the family of bounded dynamically defined operators given by~\eqref{eq.ope.din}, where $\Omega$ is a compact metric space, $T: \Omega \rightarrow \Omega$ is a minimal homeomorphism and $D,V:\Omega\rightarrow M(l,\mathbb{R})$ are continuous maps, with $D(\omega)$ invertible for each $\omega\in\Omega$. Then, 
%for each $\omega \in \Omega$, 
%$$
%\rho(H_{\omega}) \supseteq \{z \in \mathbb{C}\mid (T, A_z) \in \mathcal{UG}\},
%$$
%where the associated cocycle $(T, A_z)$ is given by \eqref{eq.cociclo.din}.
%\end{7}

\

\begin{proof5}

  Set
  \[\mathcal{A}:=\{x \in \mathbb{R}\mid (T, A(x,\cdot)) \in \mathcal{UG}\}\]
  (if $z\in\C\setminus\mathbb{R}$, then $z\in\rho(H_\omega)$, given that $H_\omega$ is self-adjoint for each $\omega\in\Omega$), and let $x\in\mathcal{A}$. Fix an arbitrary $\omega\in\Omega$. Then, %there exist constants $C>0$ and $\rho<1$ so that for each $\omega\in\Omega$, there exist a splitting $E^s\oplus E^u$ such that $\dim(E^s)=\dim(E^u)=l$,
  %$$\Vert A_n(\omega)\textbf{u}\Vert\le C\rho^n\Vert \textbf{u}\Vert,\qquad \Vert A_{-n}(\omega)\textbf{v}\Vert\le C\rho^{-n}\Vert \textbf{v}\Vert$$
%for each $n\in\mathbb{N}$, each $\textbf{u}\in E^s$ and each $\textbf{v}\in E^u$.
  %given a solution $\textbf{u}$ to the eigenvalue equation of $H_\omega$ at $z$, only one of the following conditions hold: it grows exponentially at $+\infty$, or it grows exponentially at $-\infty$. It follows from the discussion presented in Remark~\ref{lxl} that % (and so, it decays exponentially at $+\infty$). By the definition of the cocycle, it follows that
  there exist exactly $l$ linearly independent unitary vectors, %solutions to the eigenvalue equation of $H_\omega$ at $x$,
 which we denote by $\{\textbf{u}_1,\ldots,\textbf{u}_l\}$ (we omit its dependence over $x$), such that %that grow exponentially at $+\infty$ with rate $\lambda>1$ (namely,
 $\Vert A_n(x,\omega)\textbf{u}_i\Vert\le C\lambda^{-n}$ for each $i\in\{1,\ldots,l\}$ and each $n\in\N$, with another $l$ linearly independent unitary vectors, which we denote by $\{\textbf{v}_1,\ldots,\textbf{v}_l\}$, satisfying %growing exponentially at $-\infty$, with the same rate $\lambda$ (
$\Vert A_{-n}(x,\omega)\textbf{v}_i\Vert\le C\lambda^{-n}$ for each $i\in\{1,\ldots,l\}$ and each $n\in\N$ (see the proof of Theorem~\ref{UG=UH}).

%  Now, by the constancy of the Wronskian, each one of the vectors $\{\textbf{u}_1,\ldots,\textbf{u}_l\}$ (respectively, $\{\textbf{v}_1,\ldots,\textbf{v}_l\}$) decay exponentially at $-\infty$ with rate $\lambda$ (respectively, at $+ \infty$); namely, apply Green Formula~\eqref{eq.green.1} to each pair $(\textbf{v}_i,\textbf{u}_i)$, $i\in\{1,\ldots,l\}$, and use the fact that $0<\inf s_{l}[D(\omega)]\le\sup s_1[D(\omega)]<\infty$ to conclude that there exist constants $c_1,c_2>0$ such that $\Vert A_{-n}(z;\omega)\textbf{u}_i\Vert\le c_1\lambda^{-n}$ and $\Vert A_n(z;\omega)\textbf{v}_i\Vert\le c_2\lambda^{-n}$ for each $i\in\{1,\ldots,l\}$ and each $n\in\N$ (recall that $D$ is continuous).

 Thus, one may define the sequence $(F^{(+)}_{n}(x,\omega))_{n}\in M(l,\mathbb{R})$ in such way that each column of $(F^{(+)}_{n}(x,\omega)\; D_{n-1}(\omega)F^{(+)}_{n-1}(x,\omega))^t$ corresponds to $A_n(x,\omega)\textbf{v}_i$, for some $i\in\{1,\ldots,n\}$. In the same fashion, one  may define the sequence $(F^{(-)}_{n}(x,\omega))_{n}\in M(l,\mathbb{R})$ in such way that each column of $(F^{(-)}_{n}(x,\omega)\; D_{n-1}(\omega)F^{(-)}_{n-1}(x,\omega))^t$ corresponds to $A_{-n}(z,\omega)\textbf{u}_i$, for some $i\in\{1,\ldots,n\}$. Naturally, $\lim_{n\to\pm\infty}\Vert F^{(\pm)}_{n}(x,\omega)\Vert=0$.

 It follows from Proposition~\ref{prop.hip.aberto} that there exists $\epsilon>0$ such that for each $z\in B(x;\epsilon)$, $(T,A(z,\cdot))\in\mathcal{UG}$. Let $z=x+iy$, with $0<y<\epsilon/2$; %, and let $\omega\in\Omega$;
 then, by the continuity of $F^{(\pm)}_{n}(z,\omega)$ with respect to $z$ (for each fixed $n\in\Z$ and $\eta\in\Omega$, with $F^{(\pm)}_{n}(z,\eta)$ defined as above; note also that for each fixed $\eta\in\Omega$, the vectors $\{\textbf{u}_1,\ldots,\textbf{u}_l\}$ and $\{\textbf{v}_1,\ldots,\textbf{v}_l\}$ depend continuously on $z\in\C$, by the continuity of $A(z,\eta)$ with respect to $z$) and since $\lim_{n\to\pm\infty}\Vert F^{(\pm)}_{n}(z,\omega)\Vert=0$, it follows from Proposition~\ref{porp.def.green} that %for each $\omega\in\Omega$,
 \begin{eqnarray*}
   \lim_{y\downarrow 0}G_\omega(x+iy)&=&\lim_{y\downarrow 0}\sum_{j=0}^1F_{j}^{(-)}(x+iy,\omega) (Q^{\ast}(x+iy,\omega))^{-1} (F_{j}^{(+)})^{\ast}(x+iy,\omega)\\
   &=&\sum_{j=0}^1F_{j}^{(-)}(x,\omega)(Q^{t}(x,\omega))^{-1} (F_{j}^{(+)})^{t}(x,\omega)\in M(l,\R),
 \end{eqnarray*}
 and so
 \[\lim_{y\downarrow 0}\Im G_\omega(x+iy)=0.\]

 By combining this result with Definition~\ref{def.sup.minimal} and Proposition~\ref{porp.sup.ac}, one concludes that %for each $\omega\in\Omega$,
 $\mu^{tr}_\omega(\mathcal{A})=0$ (since for each $\eta\in\Omega$, $\mathcal{A}\subset \R\setminus\Sigma^\eta$ and $\mu^{tr}_\eta(\R\setminus\Sigma^\eta)=0$). By evoking again Proposition~\ref{prop.hip.aberto}, it follows that for each $x\in\mathcal{A}$ there exists $\epsilon>0$ such that $(x-\epsilon,x+\epsilon)\subset\mathcal{A}$,  so %for each $\omega\in\Omega$,
 \[\mu^{tr}_\omega((x-\epsilon,x+\epsilon))=0.\]
 This shows that %for each $\omega\in\Omega$,
 $\mathcal{A}\subset\rho(H_\omega)$ (given that $\sigma(H_\omega)=\{x\in\R\mid\forall\epsilon>0,$ $\mu^{tr}_\omega((x-\epsilon,x+\epsilon))>0\}$). Since $\omega\in\Omega$ is arbitrary, we are done.
\end{proof5}

\begin{remark} It follows from Theorem~\ref{prop.espec.3} and Remark~\ref{Rconstancy} that if $(T,\Omega)$ is transitive (with $\omega_0\in\Omega$ such that $\overline{\{T^j\omega_0\mid j\in\mathbb{Z}\}}=\Omega$), then
  \[\rho(H_{\omega_0})=\{z\in\mathbb{C}\mid(T,A(z,\cdot))\in\mathcal{UG}\}.\]
This is precisely the matrix-valued version of Theorem~1 in~\cite{alkorn}.  
\end{remark}
  
\section{Acknowledgments}

We thank the anonimous referee for invaluable suggestions that improved the quality of the manuscript. F.V. was supported by CAPES (Brazilian agency). SLC thanks the partial support by Fapemig (Minas Gerais state agency; Universal Project under contract 001/17/CEX-APQ-00352-17). 

\bibliography{bibfile}{}
\bibliographystyle{acm}

\end{document}